\documentclass[11pt]{amsart}
\usepackage{amssymb}
\usepackage{amsmath}
\usepackage{fancyhdr}
\usepackage[british]{babel}
\usepackage{geometry}
\usepackage{enumitem}
\usepackage{algpseudocode}
\usepackage{dsfont}
\usepackage{centernot}
\usepackage{xstring}
\usepackage{colortbl}
\usepackage[section]{algorithm}
\usepackage{graphicx}
\usepackage[font={footnotesize}]{caption}
\usepackage[usenames,dvipsnames,table]{xcolor}
\usepackage{pstricks,tikz}
\usepackage[h]{esvect}
\usepackage[unicode,
		bookmarksopen=true,
		bookmarksopenlevel=1,
		colorlinks=true,
		linkcolor=darkblue,
        linktoc=page,
		citecolor=darkblue,
]{hyperref}

\title[]{On a conjecture of Erd\H{o}s on locally sparse \\ Steiner triple systems}
\date{\today}
\author[S.~Glock, D.~K\"uhn, A.~Lo and D.~Osthus]{Stefan Glock, Daniela K\"uhn, Allan Lo and Deryk Osthus}

\thanks{The research leading to these results was partially supported by the EPSRC, grant nos. EP/N019504/1 (D.~K\"uhn) and EP/P002420/1 (A.~Lo),
by the Royal Society and the Wolfson Foundation (D.~K\"uhn) as well as by the European Research Council
under the European Union's Seventh Framework Programme (FP/2007--2013) / ERC Grant
Agreement no. 306349 (S.~Glock and D.~Osthus).}

\newtheorem{theorem}[algorithm]{Theorem}
\newtheorem{prop}[algorithm]{Proposition}
\newtheorem{lemma}[algorithm]{Lemma}
\newtheorem{cor}[algorithm]{Corollary}
\newtheorem{fact}[algorithm]{Fact}
\newtheorem{conj}[algorithm]{Conjecture}

\theoremstyle{definition}

\newtheorem{defin}[algorithm]{Definition}

\newtheoremstyle{claimstyle}{5pt}{5pt}{\em}{5pt}{\em}{:}{5pt}{}
\theoremstyle{claimstyle}

\newtheoremstyle{stepstyle}{10pt}{5pt}{\em}{0pt}{\em}{:}{5pt}{}
\theoremstyle{stepstyle}
\newtheorem{step}{Step}

\numberwithin{equation}{section}

\definecolor{darkblue}{rgb}{0,0,0.5}

\def\noproof{{\unskip\nobreak\hfill\penalty50\hskip2em\hbox{}\nobreak\hfill%
       $\square$\parfillskip=0pt\finalhyphendemerits=0\par}\goodbreak}
\def\endproof{\noproof\bigskip}

\newdimen\margin
\def\textno#1&#2\par{
   \margin=\hsize
   \advance\margin by -4\parindent
          \setbox1=\hbox{\sl#1}
   \ifdim\wd1 < \margin
      $$\box1\eqno#2$$
   \else
      \bigbreak
      \hbox to \hsize{\indent$\vcenter{\advance\hsize by -3\parindent
      \it\noindent#1}\hfil#2$}
      \bigbreak
   \fi}

\def\proof{\removelastskip\penalty55\medskip\noindent\setcounter{claim}{0}\setcounter{step}{0}{\bf Proof. }} 

\def\lateproof#1{\removelastskip\penalty55\medskip\noindent\setcounter{claim}{0}\setcounter{step}{0}{\bf Proof of #1. }} 

\begin{document}

\def\COMMENT#1{}
\def\TASK#1{}

\def\eps{{\varepsilon}}
\newcommand{\ex}{\mathbb{E}}
\newcommand{\pr}{\mathbb{P}}
\newcommand{\cB}{\mathcal{B}}
\newcommand{\cA}{\mathcal{A}}
\newcommand{\cE}{\mathcal{E}}
\newcommand{\cS}{\mathcal{S}}
\newcommand{\cF}{\mathcal{F}}
\newcommand{\cG}{\mathcal{G}}
\newcommand{\bL}{\mathbb{L}}
\newcommand{\bF}{\mathbb{F}}
\newcommand{\bZ}{\mathbb{Z}}
\newcommand{\cH}{\mathcal{H}}
\newcommand{\cC}{\mathcal{C}}
\newcommand{\cM}{\mathcal{M}}
\newcommand{\bN}{\mathbb{N}}
\newcommand{\bR}{\mathbb{R}}
\def\O{\mathcal{O}}
\newcommand{\cP}{\mathcal{P}}
\newcommand{\cQ}{\mathcal{Q}}
\newcommand{\cR}{\mathcal{R}}
\newcommand{\cJ}{\mathcal{J}}
\newcommand{\cL}{\mathcal{L}}
\newcommand{\cK}{\mathcal{K}}
\newcommand{\cD}{\mathcal{D}}
\newcommand{\cI}{\mathcal{I}}
\newcommand{\cV}{\mathcal{V}}
\newcommand{\cT}{\mathcal{T}}
\newcommand{\cU}{\mathcal{U}}
\newcommand{\cX}{\mathcal{X}}
\newcommand{\cZ}{\mathcal{Z}}
\newcommand{\1}{{\bf 1}_{n\not\equiv \delta}}
\newcommand{\eul}{{\rm e}}
\newcommand{\Erd}{Erd\H{o}s}
\newcommand{\cupdot}{\mathbin{\mathaccent\cdot\cup}}
\newcommand{\whp}{whp }

\newcommand{\doublesquig}{%
  \mathrel{%
    \vcenter{\offinterlineskip
      \ialign{##\cr$\rightsquigarrow$\cr\noalign{\kern-1.5pt}$\rightsquigarrow$\cr}%
    }%
  }%
}

\newcommand{\defn}{\emph}

\newcommand\restrict[1]{\raisebox{-.5ex}{$|$}_{#1}}

\newcommand{\prob}[1]{\mathrm{\mathbb{P}}\left(#1\right)}
\newcommand{\cprob}[2]{\mathrm{\mathbb{P}}_{#1}\left(#2\right)}
\newcommand{\expn}[1]{\mathrm{\mathbb{E}}\left(#1\right)}
\newcommand{\cexpn}[2]{\mathrm{\mathbb{E}}_{#1}\left(#2\right)}
\def\gnp{G_{n,p}}
\def\G{\mathcal{G}}
\def\lflr{\left\lfloor}
\def\rflr{\right\rfloor}
\def\lcl{\left\lceil}
\def\rcl{\right\rceil}

\newcommand{\qbinom}[2]{\binom{#1}{#2}_{\!q}}
\newcommand{\binomdim}[2]{\binom{#1}{#2}_{\!\dim}}

\newcommand{\grass}{\mathrm{Gr}}

\newcommand{\brackets}[1]{\left(#1\right)}
\def\sm{\setminus}
\newcommand{\Set}[1]{\{#1\}}
\newcommand{\set}[2]{\{#1\,:\;#2\}}
\newcommand{\krq}[2]{K^{(#1)}_{#2}}
\newcommand{\ind}[1]{$(\ast)_{#1}$}
\newcommand{\indcov}[1]{$(\#)_{#1}$}
\def\In{\subseteq}

\begin{abstract}  \noindent
A famous theorem of Kirkman says that there exists a Steiner triple system of order~$n$ if and only if $n\equiv 1,3\mod{6}$. In 1973, Erd\H{o}s conjectured that one can find so-called `sparse' Steiner triple systems. Roughly speaking, the aim is to have at most $j-3$ triples on every set of $j$ points, which would be best possible. 
(Triple systems with this sparseness property are also referred to as having high girth.)
We prove this conjecture asymptotically by analysing a natural generalization of the triangle removal process. Our result also solves a problem posed by Lefmann, Phelps and R\"odl 
as well as Ellis and Linial in a strong form, and answers a question of Krivelevich, Kwan, Loh and Sudakov.
Moreover, we pose a conjecture which would generalize the Erd\H{o}s conjecture to 
Steiner systems with arbitrary parameters and provide some evidence for this.
\end{abstract}

\maketitle

\section{Introduction}

Given a set $X$ of size $n$, a set $\cS$ of $3$-subsets of $X$ is a \defn{Steiner triple system of order $n$} if every $2$-subset of $X$ is contained in exactly  one of the triples of~$\cS$
(if every $2$-subset of $X$ lies in at most one of the triples of~$\cS$, we refer to $\cS$ as a \defn{partial Steiner triple 
system}).
In 1847, Kirkman~\cite{kirkman:47} proved that there exists a Steiner triple system of order~$n$ if and only if $n\equiv 1,3\mod{6}$. We shall call such $n$ \defn{admissible}.
In this paper, we investigate so-called `sparse' Steiner triple systems, which do not contain certain `forbidden configurations'. Erd\H{o}s conjectured the existence of such sparse systems. A \defn{$(j,\ell)$-configuration} is a set of $\ell$ triples on $j$ points every two of which intersect in at most one point. The `forbidden configurations' are the $(j,j-2)$-configurations. For instance, the unique $(6,4)$-configuration is called the \defn{Pasch configuration} or \defn{quadrilateral}. There are two $(7,5)$-configurations, called \defn{mitre} and \defn{mia} (see Figure~\ref{fig:examples}).

A Steiner triple system is called \defn{$k$-sparse} if it does not contain any $(j+2,j)$-configuration for $2\le j\le k$. Erd\H{o}s conjectured that if $k$ is bounded, then all these configurations can be avoided.\COMMENT{Ellis and Linial put it as: Erd\H{o}s conjectured that there are Steiner triple systems of arbitrarily high $(-2)$-girth.}

\begin{conj}[\Erd~\cite{erdos:73,erdos:76}] \label{conj:Erdos}
For every $k$, there exists an $n_k$ such that for all admissible $n>n_k$, there exists a $k$-sparse Steiner triple system of order~$n$.
\end{conj}

We note that Conjecture~\ref{conj:Erdos} would be best possible in the following sense: it is easy to see that for all $n\ge j \ge 4$, every Steiner triple system of order $n$ contains a $(j,j-3)$-configuration. This is true in a very robust sense. For instance, the $(6,3)$-theorem of Ruzsa and Szemer\'edi~\cite{RS:78} implies that any partial Steiner triple system of order $n$ with no $(6,3)$-configuration has only $o(n^2)$ triples.

The conjecture is trivial for $k\le 3$ (in the sense that it follows directly from Kirkman's theorem). A Steiner triple system is $4$-sparse if and only if it is Pasch-free. This case has received a lot of attention and has been settled in a series of papers~\cite{brouwer:77,GGW:00,GMP:90,LCGG:00}. For $k=5$, it was shown in~\cite{wolfe:05} that $5$-sparse Steiner triple systems exist for almost all admissible orders.\COMMENT{i.e.~the set of such orders has arithmetic density $1$, which is much stronger than just infinitely many}
$6$-sparse Steiner triple systems for infinitely many orders have been constructed in~\cite{FGG:07}. However, not a single $7$-sparse system is known (on at least $9$ points).
All of these and many other related results are usually based on algebraic techniques. 

Here, we prove Conjecture~\ref{conj:Erdos} approximately by analysing a natural random process. Roughly speaking, we show that when triples are randomly chosen one by one under the condition that the set of chosen triples remains sparse, then with high probability, this process runs almost to the end, i.e.~almost as many triples are added as there are in a Steiner triple system of the same order (see Theorem~\ref{thm:process stop}). In particular, such a sparse `approximate' Steiner triple system exists (the question of their existence had also been raised by Erd\H{o}s in~\cite{erdos:73}). The same result has been announced independently by Bohman and Warnke~\cite{BW:18}.

\begin{theorem}\label{thm:approx STS}
For every fixed $k$ and $n$ tending to infinity, there exists a $k$-sparse partial Steiner triple system $\cS$ on $n$ vertices with $|\cS|=(1/6-o(1))n^2$.
\end{theorem}

This also solves a problem of Lefmann, Phelps and R\"odl~\cite{LPR:93} in a strong form. They showed that for every $k$, there exists $c_k>0$ such that for all $n$ there is a $k$-sparse partial Steiner triple $\cS$ on $n$ vertices with $|\cS|\ge c_k n^2$, where $c_k\to 0$ as $k\to \infty$. Lefmann, Phelps and R\"odl asked whether $c_k$ could be bounded away from~$0$. The same question was also raised by Ellis and Linial~\cite{EL:14}. Our Theorem~\ref{thm:approx STS} implies that $c_k\sim \frac{1}{6}$ for all~$k$.
Note that the property of being $k$-sparse is often referred to as having high girth.\COMMENT{more precisely, $(-2)$-girth at least $k+3$}
Thus our result can be interpreted as providing an asymptotically optimal density bound for 
the existence of triple systems of given girth.

It is not hard to check that, for Conjecture~\ref{conj:Erdos} to be true, we must have $k=\O(\sqrt{n_k})$. In fact, $k$ needs to be much smaller than that, as shown by the following result.

\begin{theorem}[\cite{LPR:93}] \label{thm:upper bound k}
There exists $c>0$ such that every Steiner triple system of order~$n$ contains a $(j,j-2)$-configuration for some $4\le j<c \log n/\log\log n$.
\end{theorem}
This raises the question whether it is possible to allow  $k$ to grow with $n$ in Theorem~\ref{thm:approx STS}, perhaps matching the upper bound given by Theorem~\ref{thm:upper bound k}, although it is not clear what the correct function should be.

We will view configurations and partial Steiner triple systems as (linear) $3$-graphs. It will be convenient not
to assume from the outset that the systems/configurations are linear, i.e.~that every two triples meet in at most one point. Instead, we will force this condition by forbidding the so-called \defn{diamond}, i.e.~the $3$-graph with $2$ triples on $4$ vertices. 
Thus, we define a \defn{forbidden configuration} as a $3$-graph $\cS$ with $|V(\cS)|=j$ and $|\cS|=j-2$ for some $j\ge 4$. 
An \defn{\Erd-configuration} is a forbidden configuration which does not contain any forbidden configuration as a proper subgraph. Thus, the diamond is the smallest \Erd-configuration. There are no \Erd-configurations on $5$ points. Pasch and mitre are \Erd-configurations, but the mia is not as it is not Pasch-free (cf.~Figure~\ref{fig:examples}).
Clearly, a Steiner triple system is $k$-sparse if and only if it does not contain any \Erd-configuration on at most $k+2$ points. For instance, a Steiner triple system is $5$-sparse if and only if it does not contain the Pasch or the mitre configuration.
It is not too difficult to see that an \Erd-configuration exists for every order $j\ge 6$. For example, take vertices $e,o,x_1,\dots,x_{j-2}$ and all triples $ox_\ell x_{\ell+1}$ if $\ell\le j-3$ is odd and all triples $ex_\ell x_{\ell+1}$ if $\ell\le j-3$ is even. Moreover, if $j$ is even, then also take the triple $ex_{j-2}x_{1}$, and if $j$ is odd, then include the triple $x_{j-4}x_{j-2}x_{1}$ instead.\COMMENT{Let $ U $ be a vertex subset of size $u$.
Suppose that $U$ does not contain $e$ or $o$. Then $U$ contains at most $\lfloor (u-1)/2 \rfloor +1  < u - 2$ triples as $u \ge 6$. (Indeed, there is no diamond because $j\ge 6$. Thus, if $U$ induces a forbidden configuration, we must have $u\ge 6$.)
Hence $U$ contains both $e$ and $o$. 
Claim: The vertices not in $U$ are incident to at least $j - u +1$ edges. 
Since $4 \le u < j$, there exists some $r \in [j-2]$ such that $x_{r} \in U$ and $x_{r+1} \notin U$ (modulo $j-2$). 
Define $f_{\ell} = o x_\ell x_{\ell+1}$ if $\ell\le j-3$ is odd and $f_{\ell} = ex_\ell x_{\ell+1}$ if $\ell\le j-3$ is even. Moreover, set $f_{j-2} = ex_{j-2}x_{1}$ if $j$ is even, and $f_{j-2} = x_{j-4}x_{j-2}x_{1}$ if $j$ is odd.
Note that $\{ f_{\ell} \colon x_{\ell} \notin U \} \cup \{ f_{r} \}$ are $j-u+1$ distinct triples not in~$U$, proving the claim.
Hence $U$ contains at most $j -2 - (j-u+1) =  u-3$ triples, that is, $U$ does not induce a forbidden configuration.
}

\begin{figure}[ht]
\begin{center}
    \begin{tabular}{l l l l}
										$j$	&  Name      &    Triples              		&              \\ \hline
										$4^\ast$ &  diamond   &  $012,013$             &   \Erd       \\
                    $6$	&  Pasch     &  $012, 034, 135, 245$    	&    \Erd 						\\
										$7$	&  mitre     &   $012, 034, 135, 236, 456$                    		&  		\Erd					\\
										$7$	&  mia       &  $012, 034, 135, 245, 056$                     			&    	contains Pasch						\\
										$8$	&  $6$-cycle  &   $012, 034, 135, 246, 257, 367$                    &    		\Erd				\\
										$8$	&  crown     &    $012, 034, 135, 236, 147, 567$                   		&   			\Erd				\\
										$8$	&             &   $012, 034, 135, 236, 146, 057$                   		 &   	contains Pasch							 \\
										$8$	&            &    $012, 034, 135, 236, 146, 247$                  		 &    	contains Pasch								\\
										$8$	&            &    $012, 034, 135, 236, 147, 257$                    	 &    	contains mitre								\\ \hline
		\end{tabular}
\end{center}
\caption{The smallest forbidden configurations. There are more such configurations which are not linear (i.e.~contain the diamond) and thus are omitted here. If we assume at the outset that all configurations are linear, then the Pasch configuration becomes the smallest forbidden configuration.}
\label{fig:examples}
\end{figure}

As indicated above, in order to prove Theorem~\ref{thm:approx STS}, we will consider a natural random process, which can be seen as a generalization of the triangle removal process, or alternatively as an $\cH$-free process for hypergraphs.

The triangle removal process starts with the complete graph $K_n$ and then repeatedly deletes the edges of a uniformly chosen triangle. This process terminates with a triangle-free graph, and along the way produces a partial Steiner triple system. The most natural question about this process is how long it typically runs for, or equivalently, how many edges are left when it terminates. With the motivation of determining the Ramsey number $R(3,t)$, Bollob\'as and Erd\H{o}s conjectured in 1990 that with high probability the number of edges left is of order $n^{3/2}$. 
This problem attracted much attention (see e.g.~\cite{grable:97,RT:96,spencer:95}), culminating in a result of  Bohman, Frieze and Lubetzky~\cite{BFL:15}
where the exponent was finally approximately confirmed.

We adapt the triangle removal process so that it does not just produce a partial Steiner triple system, but a $k$-sparse one. Hence, in each step we delete the edges of a uniformly chosen triangle which does not produce an \Erd-configuration of order at most $k+2$ with some of the previously chosen triangles (cf.~Algorithm~\ref{alg:general removal}). The process terminates if no such triangle is left.
The question is of course again how long the process typically runs for. It was suggested by Krivelevich, Kwan, Loh, and Sudakov~\cite{KKLS:18} that the process runs for quadratically many steps.
We prove that with high probability, the number of leftover edges is $o(n^2)$, implying Theorem~\ref{thm:approx STS}. It would be interesting to find the correct order of magnitude of the number of leftover edges. It may be possible that this number is still of order $n^{3/2}$.

We actually formulate the above process as an $\cH$-free process for hypergraphs. Let $\cH$ be the set of \Erd-configurations up to order $k+2$. The \defn{$\cH$-free process} is the random process starting with an empty $3$-graph on $n$ vertices where in each step a uniformly random hyperedge is added under the condition that no copy of a member of $\cH$ is created.
For a fixed (hyper-)graph $H$, the $H$-free process has been extensively studied, in particular if $H$ is `strictly $2$-balanced' (see e.g.~\cite{BB:16,BK:10,KOT:16,OT:01,warnke:14b,warnke:14a}).
A particular challenge arising in the analysis of the current process
is that each individual Erd\H{o}s-configuration in $\cH$ has a significant influence on the trajectory of the process.

An obvious question is whether our approximate result can be combined with the absorbing method in order to prove Conjecture~\ref{conj:Erdos}, e.g.~using approaches from \cite{keevash:14,keevash:18} or~\cite{GKLO:16}. One major difficulty here is that the absorbing method relies on the simple fact that, given two triangle packings which are edge-disjoint, their union also forms a triangle packing. On the contrary, the union of two edge-disjoint sparse triangle packings is not necessarily sparse.

Our paper is organised as follows. After introducing our basic terminology in Section~\ref{sec:notation}, we will state Freedman's inequality in Section~\ref{sec:freedman}, which will be the main probabilistic tool to analyse our process. In Section~\ref{sec:process}, we define the process more formally, discuss the key random variables of the process and predict its behaviour heuristically using the differential equation method. Subsequently, in Section~\ref{sec:analysis}, we analyse the process. In particular, we establish trend hypotheses and boundedness hypotheses for the random variables which we track.
In Section~\ref{sec:counting}, we formulate a conjecture on the number of $k$-sparse Steiner triple systems. Finally, in Section~\ref{sec:designs}, we propose a conjecture which would generalize Conjecture~\ref{conj:Erdos} to Steiner systems with arbitrary parameters and provide some evidence for our conjecture.

\section{Notation} \label{sec:notation}

We let $[n]$ denote the set $\Set{1,\dots,n}$, where $[0]:=\emptyset$. Moreover, we set $[n]_0:=[n]\cup\Set{0}$ and $\bN_0:=\bN\cup \Set{0}$. Given a set $X$ and $i\in\bN_0$, we write $\binom{X}{i}$ for the collection of all $i$-subsets of $X$.

A \defn{hypergraph} $H$ is a pair $(V,E)$, where $V=V(H)$ is the vertex set and the edge set $E$ is a set of subsets of $V$. We identify $H$ with $E$.
In particular, we let $|H|:=|E|$. We say that $H$ is an \defn{$r$-graph} if every edge has size $r$. Given $U\subseteq V(H)$, we write $H[U]$ for the sub-hypergraph of $H$ induced by $U$. Given $S\subseteq V(H)$, we write $d_H(S)$ for the \defn{degree of $S$ in $H$}, i.e.~the number of hyperedges of $H$ containing $S$.

We say that an event holds \defn{with high probability (whp)} if the probability that it holds tends to $1$ as $n\to\infty$ (where $n$ usually denotes the number of vertices). 

We write $a=b\pm c$ if $b-c\le a\le b+c$. Equations containing $\pm$ are always to be interpreted from left to right, e.g. $b_1\pm c_1=b_2\pm c_2$ means that $b_1-c_1\ge b_2-c_2$ and $b_1+c_1\le b_2+c_2$.
Moreover, $a\wedge b$ denotes the minimum of $a$ and~$b$.

We write $f=\O(g)$ if $|f|\le C|g|$ for some constant $C$ (which by default may only depend on $k$).\COMMENT{If $f$ and $g$ are functions of $x$, we write $f(x)=\O(g(x))$ for $x\in D$ if $|f(x)|\le C|g(x)|$ for all $x\in D$.}  We write $\O_{\gamma}$ to indicate that the constant may also depend on~$\gamma$. Similarly, we write $f=\Omega(g)$ if $f\ge c|g|$ for some constant $c>0$ (which by default may depend only on $k$, and additional dependencies are indicated as indices).
Note that if $x=\O(n^a)$ and $y=\Omega(n^b)$, then
\begin{align}
\frac{x+ \O(\eps) n^a}{y+ \O(\eps) n^b}=\frac{x}{y}+ \O(\eps)n^{a-b}.\label{O fractions}
\end{align}

We write $x\ll y$ to mean that for any $y\in (0,1]$ there exists an $x_0\in (0,1)$ such that for all $x\le x_0$ the subsequent statement holds. Hierarchies with more constants are defined in a similar way and are to be read from the right to the left. We will always assume that the constants in our hierarchies are reals in $(0,1]$. Moreover, if $1/x$ appears in a hierarchy, this implicitly means that $x$ is a natural number. More precisely, $1/x\ll y$ means that for any $y\in (0,1]$ there exists an $x_0\in \bN$ such that for all $x\in \bN$ with $x\ge x_0$ the subsequent statement holds.

\section{Freedman's inequality} \label{sec:freedman}

Let $X(0),X(1),\dots$ be a real-valued random process. We define $$\Delta X(i):=X(i+1)-X(i).$$ The process $X(0),X(1),\dots$ is a \defn{supermartingale} (with respect to a filtration $\cF=(\cF(0),\cF(1),\dots)$) if $\expn{X(i+1)\mid \cF(i)}\le X(i)$, or equivalently, $\expn{\Delta X(i)\mid \cF(i)}\le 0$, for all~$i\ge 0$.

The following tail probability is due to Freedman~\cite{freedman:75}. It was originally stated for martingales, but the proof for supermartingales is verbatim the same.
\begin{lemma}[Freedman's inequality~\cite{freedman:75}] \label{lem:freedman}
Let $X(0),X(1),\dots$ be a supermartingale with respect to a filtration $\cF=(\cF(0),\cF(1),\dots)$.
Suppose that $\left|\Delta X(i)\right|\le K$ for all $i$, and let $V(i):=\sum_{j=0}^{i-1}\expn{\left(\Delta X(j)\right)^{2}\mid \cF(j)}$.
Then for any $t,v>0$,
$$\prob{X(i)\ge X(0)+t  \mbox{ and } V(i)\le v \mbox{ for some }i}\le \eul^{-\frac{t^{2}}{2\left(v+Kt\right)}}.$$
\end{lemma}

We will apply Lemma~\ref{lem:freedman} in the following scenario: There will be a parameter $n$ which measures the size of the probability space. There will be a (random) time $\tau_{freeze}=\O(n^2)$ such that $\Delta X(i)=0$ for all $i\ge \tau_{freeze}$. Moreover, we will have $\Delta X(i)= \O(n^{\alpha_2})$ and $\expn{|\Delta X(i)| \mid \cF(i)}=\O(n^{\alpha_3})$ for all $i$, and $-X(0)=\Omega(n^{\alpha_1})$. Suppose that $\alpha_1>\alpha_2$ and $\alpha_1\ge \alpha_3+2$. Then we can conclude that
\begin{align}
\prob{\exists i\colon X(i)\ge 0}\le \eul^{-\Omega(n^{\alpha_1-\alpha_2})}.\label{freedman compact}
\end{align}
Indeed, we can apply Lemma~\ref{lem:freedman} with $t:=-X(0)$, $v=\O(n^{\alpha_2+\alpha_3+2})$ and $K=\O(n^{\alpha_2})$. For every $i$, we have $\expn{\left(\Delta X(i)\right)^{2}\mid \cF(i)}\le K \cdot  \expn{\left|\Delta X(i)\right|\mid \cF(i)}=\O(n^{\alpha_2+\alpha_3})$ and thus $\sum_{i=0}^{\infty}\expn{\left(\Delta X(i)\right)^{2}\mid \cF(i)}\le v$. Hence, $V(i)\le v$ for all~$i$. Note that $t^2=\Omega(n^{2\alpha_1})$ and $v+Kt=\O(n^{\alpha_1+\alpha_2})$.

The supermartingales we consider are obtained as follows: Let $X$ be a random variable of the process, e.g.~the number of available triples containing a fixed edge. Using the differential equation method, we will have a rough idea of how $X$ should behave, i.e.~we will find a (smooth) deterministic function $f_X$ and predict that $X\approx f_X$. We call $f_X$ the \defn{trajectory} of $X$. In order to control the deviation of $X$ from $f_X$, we introduce an \defn{error function} $\eps_X$. We now define
\begin{align*}
X^+(i) &:=X(i)-f_X(i)-\eps_X(i),\\
X^-(i) &:=-X(i)+f_X(i)-\eps_X(i).
\end{align*}
(In the actual proof we will actually `freeze' these variables after a certain random time.)
Note that if $X^{\pm}(i)\le 0$, then $|X(i)-f_X(i)|\le \eps_X(i)$. (We write $X^{\pm}(i)\le 0$ to mean that both $X^{+}(i)\le 0$ and $X^{-}(i)\le 0$ hold.) Our aim is to show that the $X^{\pm}$ define two supermartingales and then to use~\eqref{freedman compact} to show that $X^{\pm}\le 0$ throughout the process.
In order to show that $X^{\pm}$ are supermartingales (with respect to a filtration $\cF=(\cF(0),\cF(1),\dots)$), it is enough to show that $\expn{\Delta X^{\pm}(i)\mid \cF(i)}\le 0$ for all $i\ge 0$ (usually referred to as the `trend hypothesis').
Observe that $$\expn{\Delta X^{\pm}(i)\mid \cF(i)}= \pm \expn{\Delta X(i)\mid \cF(i)} \mp \Delta f_X(i) -\Delta \eps_X(i).$$
In order to determine $\Delta f_X$ and $\Delta \eps_X$, we use the following simple consequence of Taylor's theorem with remainder in Lagrange form: for a sufficiently smooth function~$f$, we have
\begin{align}
\Delta f(i) := f(i+1)-f(i)=f'(i)\pm \sup_{\xi\in[i,i+1]}f''(\xi). \label{Taylor}
\end{align}
The terms $\expn{\Delta X(i)\mid \cF(i)}$ and $\Delta f_X(i)$ will almost cancel out, and the purpose of $\Delta \eps_X(i)$ is to make the sum negative. For this to work, $\eps_X$ has to have a large enough growth rate throughout the process. On the other hand, it must not grow too fast, otherwise we would lose control of $X$. A careful calibration is thus essential for the analysis to work.

Once we have established that $X^+$ is a supermartingale, it remains to give bounds on $|\Delta X^+(i)|$ (`boundedness hypothesis') and $\expn{|\Delta X^+(i)| \mid \cF(i)}$. For this, we simply use $|\Delta X^+(i)|\le |\Delta X(i)|+|\Delta f_X(i)|+|\Delta \eps_X(i)|$.

\COMMENT{Could delete the rest}Let $\cX(i)$ be a set (which contains all objects of a certain type at time $i$) and suppose that our random variable is defined as $X(i):=|\cX(i)|$.
Suppose we consider our process at time $i$. For every object $x\in \cX(i)$, $x$ could potentially be removed from $\cX(i)$, i.e.~$x\notin \cX(i+1)$. We denote the indicator function of this event by $\mathds{1}_{-x}$. Moreover, there is a set $\cX^{pot}(i)$ of potential new elements which might be added to $\cX(i)$, i.e.~for every $x\in\cX^{pot}(i)$, we have $x\notin \cX(i)$ but with non-zero probability we have $x\in \cX(i+1)$. We denote the indicator function of this event by $\mathds{1}_{+x}$. Thus, we have
\begin{align}
\Delta X(i)=|\cX(i+1)|-|\cX(i)|= - \sum_{x\in \cX(i)}\mathds{1}_{-x}   +    \sum_{x\in \cX^{pot}(i)}\mathds{1}_{+x}.
\end{align}

\section{The process} \label{sec:process}

We now describe the process that we wish to analyse. Let $V$ be a set of $n$ vertices. Suppose that we want to construct a $k$-sparse triple system, with $k\ge 2$.
Let $$j_{max}:=k+2$$
and consider Algorithm~\ref{alg:general removal}.
\begin{algorithm}
\caption{}
\label{alg:general removal}
\begin{algorithmic}
\State{$\cA(0):=\binom{V}{3}$, $\cC(0):=\emptyset$, $i:=0$}
\While{$\cA(i)\neq \emptyset$}
\State{select $T^\ast(i)\in \cA(i)$ uniformly at random}
\State{let $\cA'(i)$ consist of all $T\in \cA(i)$ for which there is $\cC' \In \cC(i)$ such that $\Set{T,T^\ast(i)}\cup \cC'$ is an \Erd-configuration on at most $j_{max}$ points}
\State{$\cA(i+1):=\cA(i)\sm (\cA'(i)\cup \Set{T^\ast(i)})$}
\State{$\cC(i+1):=\cC(i)\cup \Set{T^\ast(i)}$}
\State{$i:=i+1$}
\EndWhile
\end{algorithmic}
\end{algorithm}

Clearly, it is enough to forbid \Erd-configurations on at most $j_{max}$ points. In our analysis, it will be important that we only consider these `minimal' forbidden configurations, as it turns out that they behave `almost independently', which would not be the case if we considered all forbidden configurations.

The last step of the process is $\tau_{max}:=\min\set{i}{\cA(i)=\emptyset}$.
At time $i$, we say that $\cA(i)$ is the set of \defn{available triples} and $\cC(i)$ is the set of \defn{chosen triples}. Clearly, we have $\cA(i+1)\In \cA(i)$, $\cC(i+1)\supseteq \cC(i)$ and $\cA(i)\cap \cC(i)=\emptyset$ for all $i$. We refer to $T^\ast(i)$ as the \defn{selected triple} in step $i$.
For a $3$-set $T\In V$, let $\tau_T:=\min\set{i}{T\notin \cA(i)}$.

\begin{fact} \label{fact:process works}
$|\cC(i)|=i$ and $\cC(i)$ is $k$-sparse for all $i\le \tau_{max}$.
\end{fact}

In particular, $\cC(i)$ is a linear $3$-graph, i.e.~$|T^\ast(i')\cap T^\ast(i'')|\le 1$ for all distinct $i',i''<i$.

A $2$-set $e\In V$ is called \defn{covered} (at time $i$) if $e\In T$ for some $T\in \cC(i)$,\COMMENT{before: there is $i'<i$ such that $e\In T^\ast(i')$} otherwise it is \defn{uncovered}. We often refer to $2$-sets of $V$ as \defn{edges}. Let $E(i)$ be the set of uncovered edges at time $i$. Since $\cC(i)$ is linear, we have $|E(i)|=\binom{n}{2}-3|\cC(i)|=\binom{n}{2}-3i$ for all $i\le \tau_{max}$.
For a $2$-set $e$, we define the random time $\tau_e:=\min\set{i}{e\notin E(i)}$, where $\tau_e:=\infty$ if $e\in E(\tau_{max})$.
\begin{fact}\label{fact:edge covered}
If $T\in \binom{V}{3}$ is available, then every edge contained in $T$ is uncovered.
\end{fact}

By Fact~\ref{fact:process works}, the following result implies Theorem~\ref{thm:approx STS}.

\begin{theorem}\label{thm:process stop}
Suppose that $\gamma\in(0,1)$ and $k\in \bN$. Then whp as $n\to \infty$, $\tau_{max} \ge (1-\gamma)n^2/6$.
\end{theorem}

\subsection{Key variables and threats} \label{subsec:key vars}

Define the densities
\begin{align}
p(i) :=|E(i)|/\binom{n}{2},\quad p_{\cC}(i) :=|\cC(i)|/\binom{n}{3},\quad p_{\cA}(i) &:=|\cA(i)|/\binom{n}{3}. \label{def densities}
\end{align}

The following equalities clearly hold throughout the process, i.e.~for all $i\le \tau_{max}$:
\begin{align}
p(i) & =1-\frac{3i}{\binom{n}{2}},  \label{edge density}\\ 
p_{\cC}(i) &=\frac{i}{\binom{n}{3}}.
\end{align}

However, this gives no information as to how long the process continues. For this, we need to track the number $|\cA(i)|$ of available triples.

For $T_1,T_2\in \cA(i)$, we say that $T_1$ and $T_2$ \defn{exclude each other}, denoted by $T_1 \leftrightarrow T_2$, if there is $\cC' \In \cC(i)$ such that $\Set{T_1,T_2}\cup \cC'$ is an \Erd-configuration on at most $j_{max}$ points.

For $T\in \cA(i)$, let $\cT_T(i):=\set{T^\ast\in \cA(i)}{T \leftrightarrow T^\ast}$. 
Hence, if $T^\ast(i)\in \cT_T(i)$ then $T\in \cA'(i)$. Note that $T\notin \cT_T(i)$.\COMMENT{Otherwise wouldn't be available anymore}
Since $T^\ast(i)$ is selected uniformly at random, we have that the probability that $T$ is not in $\cA(i+1)$ is $\frac{|\cT_T(i)|+1}{|\cA(i)|}$.

For a $2$-set~$e$, let $\cX_e(i):=\set{T\in \cA(i)}{T\supseteq e}$ be the set of available triples containing $e$ at time~$i$.
Moreover, we set $X_e(i):=|\cX_e(i)|$. Clearly, we have $X_e(0)=n-2$.

\begin{fact} \label{fact:available edge triple relation}
$|\cA(i)|=\frac{1}{3}\sum_{e\in E(i)} X_e(i)$.
\end{fact}

\proof
By Fact~\ref{fact:edge covered}, every available triple contains $3$ uncovered edges, and $X_e(i)=|\cX_e(i)|$ for all $e\in E(i)$.
\endproof

Let $\mathfrak{J}_j$ be the set of all unlabelled\COMMENT{i.e.~only care about the set of triples} \Erd-configurations on $j$ vertices in $V$. For a triple~$T$, we let $\mathfrak{J}_j(T):=\set{\cS\in \mathfrak{J}_j}{T\in \cS}$. By symmetry, we have that $|\mathfrak{J}_j(T)|=:J_j$ is the same for all triples $T$. We will not compute the precise number, but only need that $J_j=\Theta(n^{j-3})$ for $j\ge 6$.\COMMENT{can write $J_j=erd_j\binom{n-3}{j-3}$, see Section~\ref{sec:counting}.}

For a triple $T$, $j\in\Set{4,\dots,j_{max}}$ and $c\in\Set{0,\dots,j-4}$, we define
\begin{align}
\cX_{T,j,c}(i):=\set{\cS\in \mathfrak{J}_j(T)}{|(\cS-\Set{T})\cap \cC(i)|=c,|(\cS-\Set{T})\cap \cA(i)|=j-3-c}. \label{def erdos}
\end{align}
Note that if $\cS\in \cX_{T,j,c}(i)$, then every $T'\in \cS-\Set{T}$ is either chosen or available (at time $i$). We make no assumption on the status of $T$, however we will only be interested in $\cX_{T,j,c}(i)$ as long as $T$ is available.
Define $X_{T,j,c}(i):=|\cX_{T,j,c}(i)|$. Note that $X_{T,j,0}(0)=J_j$ and $X_{T,j,c}(0)=0$ if $c>0$.

Note that since $\mathfrak{J}_5=\emptyset$, we always have $\cX_{T,5,c}(i)=\emptyset$. Moreover, note that $\cX_{T,4,0}(i)$ corresponds to the set of all $T'\in \cA(i)$ with $|T'\cap T|=2$.

We call elements of $\cX_{T,j,j-4}$ \defn{dangerous configurations}.
\begin{fact} \label{fact:threats via extensions}
For $T\in \cA(i)$, we have $$\cT_T(i)=\set{T^\ast}{\exists\, \cS\in \bigcup_{j=4}^{j_{max}}\cX_{T,j,j-4}(i) \mbox{ such that } (\cS-\Set{T})\cap \cA(i)=\Set{T^\ast}}.$$
\end{fact}

\proof
Suppose $T^\ast\in \cT_T(i)\In \cA(i)$. Then there is $\cC'\In \cC(i)$ such that $\Set{T,T^\ast}\cup \cC'$ forms an \Erd-configuration $\cS$ on $j\le j_{max}$ points. Then $\cS\in \cX_{T,j,j-4}(i)$ and $(\cS-\Set{T})\cap \cA(i)=\Set{T^\ast}$. Conversely, if there is $\cS\in \bigcup_{j=4}^{j_{max}}\cX_{T,j,j-4}(i)$ with $(\cS-\Set{T})\cap \cA(i)=\Set{T^\ast}$, then $\cC':=\cS-\Set{T,T^\ast}\In \cC(i)$ is such that $\Set{T,T^\ast}\cup \cC'$ forms an \Erd-configuration on at most $j_{max}$ points.
\endproof

By showing that most $T^\ast$ are only contained in at most one $\cS\in \bigcup_{j=4}^{j_{max}}\cX_{T,j,j-4}(i)$, we will see (cf.~Proposition~\ref{prop:danger count}) that
\begin{align}
|\cT_T(i)| &\approx \sum_{j=4}^{j_{max}} X_{T,j,j-4}(i).\label{threat approximation}
\end{align}

\begin{fact}\label{fact:smallest Erdos}
For $T\in \cA(i)$, $X_{T,4,0}(i) = \sum_{e\in \binom{T}{2}}X_e(i)-3$.
\end{fact}

\proof
For $T\in \cA(i)$, $X_{T,4,0}(i)$ counts the number of $T^\ast\in \cA(i)$ with $|T\cap T^\ast|=2$. If for such $T^\ast$, we have $T\cap T^\ast=e \in \binom{T}{2}$, then $T^\ast \in \cX_e(i)\sm \Set{T}$. 
\endproof

For $e\in E(i)$ and $T\in \cX_e(i)$, we say that $T^\ast\in \cA(i)$ \defn{threatens~$T,e$} if $e\not\In T^\ast$ and $T \leftrightarrow T^\ast$. This means that if $T^\ast$ is the selected triple $T^\ast(i)$, then $T\notin \cA(i+1)$, but still $e\in E(i+1)$. Let $th_{T,e}(i)$ be the number of threats to $T,e$.

\begin{prop}\label{prop:edge triple dangers}
For $e\in E(i)$ and $T\in \cX_e(i)$, we have $th_{T,e}(i)=|\cT_T(i)|-X_e(i)+1$.
\end{prop}

\proof
We have $th_{T,e}(i)=|\cT_T(i)|-|\set{T^\ast\in \cT_T(i)}{e\In T^\ast}|$. Since for $j\ge 5$ and $T^\ast\in \cT_T(i)$ with $e\In T^\ast$, there is no \Erd-configuration on $j$ points which contains $T$ and $T^\ast$,\COMMENT{as $\Set{T,T^\ast}$ forms a diamond} we have $\set{T^\ast\in \cT_T(i)}{e\In T^\ast}=\set{T^\ast\in \cA(i)\sm\Set{T}}{e\In T^\ast}=\cX_e(i)\sm\Set{T}$.
\endproof

Together with~\eqref{threat approximation} and Fact~\ref{fact:smallest Erdos}, we have that 
\begin{align}
th_{T,e}(i) &\approx \sum_{e'\in \binom{T}{2}\sm\Set{e}}X_{e'}(i)+ \sum_{j=6}^{j_{max}} X_{T,j,j-4}(i).\label{threat approximation edge}
\end{align}

For $T\in \cA(i)$ and $\cS\in \cX_{T,j,c}(i)$, we say that $T^\ast\in \cA(i)$ \defn{threatens} $\cS,T$ if $T \centernot\leftrightarrow T^\ast $ and $T\neq T^\ast$ and there is $T'\in (\cS-\Set{T})\cap \cA(i)$ with $T'\leftrightarrow T^\ast$ or $T'=T^\ast$. (Note that if $c=j-4$, then the case $T'=T^\ast$ cannot happen as this would imply $T\leftrightarrow T^\ast $.) This means that if $T^\ast$ is the selected triple $T^\ast(i)$, then $\cS\notin \cX_{T,j,c}(i+1)$, but still $T\in \cA(i+1)$. Let $th_{\cS,T}(i)$ be the number of threats to $\cS,T$.

By showing that, for fixed $\cS,T$, most $T^\ast$ exclude only one $T'\in (\cS-\{T\})\cap \cA(i)$, we will see (cf.~Proposition~\ref{prop:Erdos threats}) that
\begin{align}
th_{\cS,T}(i) &\approx \sum_{T'\in (\cS-\Set{T})\cap \cA(i)}|\cT_{T'}(i)|.\label{threat approximation triple}
\end{align}

\subsection{Heuristics} \label{subsec:heuristics}

We now heuristically predict the behaviour of the process. This is only a heuristic argument and not a part of the formal proof, yet should provide motivation for our choice of the trajectories $f_X$. As part of the exposition, we define some key functions which will play a crucial role in the remainder of the paper.

We make the assumptions that for all $e\in E(i)$, we have $X_e(i)\approx f_{edge}(i)$
for some function $f_{edge}$. Similarly, for all $T\in \cA(i)$ and $j\in\Set{6,\dots,j_{max}}$, we have $X_{T,j,j-4}(i)\approx f_{j,j-4}(i)$ for some function $f_{j,j-4}$. In order to use the differential equation method, we interpret the term $\dfrac{\mathrm{d}f_{edge}(i)}{\mathrm{d}i}$ as the expectation of $\Delta X_e(i)$. Note that $f_{edge}(i)$ approximates $X_e(i)$ only for uncovered edges~$e$, whereas for all covered edges $e$ we have $X_e(i)=0$. Thus, it is important to consider the conditional expectation of $\Delta X_e(i)$ under the event that $e$ remains uncovered. In this conditional probability space, $T^\ast(i)$ is chosen uniformly from $\cA(i)\sm \cX_e(i)$, and for a fixed triple $T\in \cX_{e}(i)$, we have $T\notin \cX_{e}(i+1)$ if and only if $T^\ast(i)$ threatens $T,e$. Thus, the (conditional) probability that $T$ becomes unavailable is given by $\frac{th_{T,e}(i)}{|\cA(i)\sm \cX_e(i)|}$.
For brevity, define 
\begin{align}
A(i) &=\frac{1}{3}|E(i)|f_{edge}(i) \overset{\eqref{def densities}}{=} \frac{1}{3}p(i)\binom{n}{2}f_{edge}(i), \label{available triple edge rel heur} \\ 
F(i) &=\sum_{j=6}^{j_{max}} f_{j,j-4}(i). \label{def sum F heur} 
\end{align}
Fact~\ref{fact:available edge triple relation} indicates that $|\cA(i)\sm \cX_e(i)| \approx |\cA(i)|\approx A(i)$. From~\eqref{threat approximation edge}, we deduce that $th_{T,e}(i)\approx 2f_{edge}(i)+F(i)$.
Thus, we approximate the conditional expectation of $\Delta X_e(i)$ as $$-\sum_{T\in \cX_e(i)}\frac{th_{T,e}(i)}{|\cA(i)\sm \cX_e(i)|}\approx - \frac{2f_{edge}(i)+F(i)}{A(i)}f_{edge}(i). $$
We obtain the following differential equation for $f_{edge}(i)$:
\begin{align}
\dfrac{\mathrm{d}f_{edge}(i)}{\mathrm{d}i}=- \frac{2f_{edge}(i)+F(i)}{A(i)}f_{edge}(i). \label{diff eq edges}
\end{align}

In order to obtain an expression for $F(i)$, we make the additional assumption that $\cA(i)$ and $\cC(i)$ are random $3$-graphs obtained by including every triple independently with probability $p_{\cA}$ and $p_{\cC}$, respectively, conditioned on $\cA(i)\cap \cC(i)=\emptyset$. Fix a triple $T\In V$. Recall that $J_j=|\mathfrak{J}_j(T)|$ denotes the number of unlabelled \Erd-configurations on $j$ points in $V$ which contain $T$ as a triple. For each $\cS\in \mathfrak{J}_j(T)$, $\cS$ belongs to $\cX_{T,j,j-4}$ if and only if one triple of $\cS-\Set{T}$ is available, and the other $j-4$ are chosen. Under the above assumption, the probability for this is $(j-3) p_{\cC}(i)^{j-4} p_{\cA}(i)$.
We thus guess that
\begin{align}
f_{j,j-4}(i)= (j-3)p_{\cC}(i)^{j-4} p_{\cA}(i) J_j.\label{naive XT prediction}
\end{align}
By \eqref{def densities} and \eqref{def sum F heur}, we then have $$F(i)=\sum_{j=6}^{j_{max}}(j-3)\left(\frac{i}{\binom{n}{3}}\right)^{j-4} \frac{|\cA(i)|}{\binom{n}{3}}J_j=|\cA(i)|\sum_{j=6}^{j_{max}}\frac{(j-3)J_j}{\binom{n}{3}^{j-3}}i^{j-4}.$$
This motivates the definition of the following function, which turns out to be a crucial parameter of the process.
\begin{align*}
\rho'(i)=\sum_{j=6}^{j_{max}}\frac{(j-3)J_j}{\binom{n}{3}^{j-3}}i^{j-4}. 
\end{align*}
We obtain
$$\frac{2f_{edge}(i)+F(i)}{A(i)} \approx \frac{2f_{edge}(i)+A(i)\rho'(i)}{A(i)} = \frac{6}{p(i)\binom{n}{2}}+\rho'(i).$$
Substituting this into \eqref{diff eq edges} yields the linear differential equation
\begin{align*}
\dfrac{\mathrm{d}f_{edge}(i)}{\mathrm{d}i}=-f_{edge}(i)\left(\frac{6}{p(i)\binom{n}{2}}+\rho'(i)\right).
\end{align*}
For this equation, we can find the solution (e.g.~using separation of variables)
\begin{align*}
f_{edge}(i)=\eul^{-\rho(i)}p(i)^2 f_{edge}(0)=\eul^{-\rho(i)}p(i)^2 (n-2), 
\end{align*}
where
\begin{align*}
\rho(i)=\sum_{j=6}^{j_{max}}\frac{J_j}{\binom{n}{3}^{j-3}}i^{j-3}
\end{align*}
is the integral of $\rho'$ with $\rho(0)=0$.

We briefly interpret this result. Note that since $J_j=\Theta(n^{j-3})$ and $i=\O(n^2)$, we have that $\rho(i)=\O(1)$. Also, as long as $i=o(n^2)$, we have $\rho(i)=o(1)$, i.e.~the effect of the term $\eul^{-\rho(i)}$ is negligible. This means that in the early stages of the process, we expect $X_e$ to behave as in the standard random triangle removal process. Once $i$ is quadratic in~$n$, sufficiently many dangerous configurations have been created to affect $X_e(i)$ significantly. However, their influence is limited in the sense that they modify $X_e(i)$ only by a multiplicative constant.

\subsection{Trajectories}

As a result of the heuristic argument, we conclude that we wish to track the random variables $X_e$ and $X_{T,j,j-4}$. Clearly, in order to track $X_{T,j,j-4}$, we also need to track $X_{T,j,j-5}$, and so on. A guess for the trajectory of $X_{T,j,c}$ can be obtained similarly to \eqref{naive XT prediction}.
We now define the trajectories for these key variables formally. For clarity, we also define the other relevant functions from above again.

\begin{defin}[Trajectories]
For $0\le i\le \binom{n}{2}/3$, define the functions
\begin{align}
\rho(i)      &:=\sum_{j=6}^{j_{max}}\frac{J_j}{\binom{n}{3}^{j-3}}i^{j-3},    \label{def rho} \\
f_{edge}(i) 			&:= \eul^{-\rho(i)}p(i)^2 (n-2),       \label{edge trajectory}  \\
A(i)       &:=\eul^{-\rho(i)}p(i)^3 \binom{n}{3} \quad \left(=\frac{1}{3}p(i)\binom{n}{2}f_{edge}(i)\right).    \label{available trajectory}
\end{align}
Moreover, for all $j\in \Set{6,\dots,j_{max}}$, define
\begin{align}
f_{j,c}(i) &:= \binom{j-3}{c}\eul^{-(j-3-c)\rho(i)}p(i)^{3(j-3-c)}i^c \binom{n}{3}^{-c}J_j \label{triple trajectory}
\end{align}
for all $c\in\Set{1,\dots,j-4}$ and $f_{j,0}(i):= \eul^{-(j-3)\rho(i)}p(i)^{3(j-3)}J_j$.\COMMENT{bypasses $0^0$ if $i=0$ and makes it clearer when computing derivatives}
Finally, define
\begin{align}
F(i)        &:=\sum_{j=6}^{j_{max}} f_{j,j-4}(i). \label{def sum F}
\end{align}
\end{defin}

We close this section by observing the following properties of the functions defined above.

\begin{prop}
For $0\le i\le \binom{n}{2}/3$, the following hold:
\begin{align}
&\rho(i)=\O(1), \; \rho'(i)=\O(n^{-2}), \; \rho''(i)=\O(n^{-4}); \label{rho magnitudes} \\
&f_{j,c}(i)=\O(n^{j-3-c}) \mbox{ for all $j,c$};  \label{triple trajectory magnitude} \\
&f_{edge}(i)=\O(n), \; F(i)=\O(n).  \label{edge and danger magnitude}
\end{align}
\end{prop}

\proof
For \eqref{rho magnitudes}, recall that $J_j=\Theta(n^{j-3})$. \eqref{triple trajectory magnitude} and~\eqref{edge and danger magnitude} then follow.
\endproof

\section{Analysis of the process} \label{sec:analysis}

In this section, we prove Theorem~\ref{thm:process stop}.
Choose constants
\begin{align}
\eps_0 \ll 1/C \ll \gamma \ll 1/k. \label{hierarchy}
\end{align} In all calculations, we assume that $n$ is sufficiently large once all other constants are fixed.

\subsection{Extension types}
As mentioned before, in order to track our key variables, we also need to track a number of auxiliary variables, e.g.~to account for double counting when estimating the threat functions $th_{T,e}$ and $th_{\cS,T}$. We will also need such auxiliary variables to establish boundedness conditions for our key variables. Fortunately, it suffices to have (generous) upper bounds on these variables. This allows us to treat all the auxiliary variables we need using a unified framework.

An \defn{extension type} is a pair $(H,U)$ where $H$ is a $3$-graph and $U\In V(H)$ is such that $|H[U]|=0$. We can think of $U$ as a set of \defn{root vertices}, whereas the vertices in $V(H)\sm U$ are \defn{free}.
Given a $3$-graph $G$ and a set $R\In V(G)$ with $|R|=|U|$, an \defn{$(H,U)$-extension at $R$ in $G$} is an embedding $\phi\colon H \to G$ such that $\phi(U)=R$, i.e.~an injective map $\phi\colon V(H) \to V(G)$ such that $\phi(U)=R$ and $\phi(e)\in G$ for all $e\in H$. Note that if $G$ is a random $3$-graph on $n$ vertices, where edges appear independently with probability $1/n$, then the expected number of $(H,U)$-extensions at a fixed set $R$\COMMENT{of constant size} is of order $n^{|V(H)\sm U|-|H|}$.

Let $$m:=2j_{max}.$$
For an extension type $(H,U)$ with $|V(H)|\le m$ and a set $R\In V$ with $|R|=|U|$, we define the random variable $X_{R,(H,U)}(i)$ counting the number of $(H,U)$-extensions at $R$ in $\cC(i)$. Note that $X_{R,(H,U)}(0)=0$ if $H$ is non-empty.

\begin{defin}
Call an extension type $(H,U)$ \defn{$\kappa$-balanced} if for all $U\In U'\In V(H)$, we have $|H-H[U']|\ge |V(H)\sm U'|-\kappa$. Let $\kappa(H,U):=\min\set{\kappa\ge 0}{(H,U)\mbox{ is }\kappa\mbox{-balanced}}$.\COMMENT{Should we mention a connection to strictly $k$-balanced definition?}

For $\ell\in \Set{1,\dots,m}$ and $\kappa \in\Set{0,\dots,\ell}$, let $Ext(\kappa,\ell)$ denote the set of all extension types $(H,U)$ with $|V(H)|\le m$, $\kappa(H,U)=\kappa$ and $|V(H)\sm U|=\ell$, and such that $H$ is not empty. We do not distinguish between isomorphic extension types here. In particular, $|Ext(\kappa,\ell)|=\O(1)$.
\end{defin}

We gather a few easy facts about balanced extension types.
\begin{fact} \label{fact:extension type trivial}
Let $(H,U)$ be an extension type, $\kappa:=\kappa(H,U)$ and $\ell:=|V(H)\sm U|$. Then the following assertions hold.
\begin{enumerate}[label=\rm{(\roman*)}]
\item $|V(H)\sm U|-|H|\le \kappa$. \label{fact:extension type l kappa}
\item $\kappa\le \ell$.
\item If $H$ is empty, then $\kappa=\ell$. \label{fact:extension type empty}
\item If $R\In V$ with $|R|=|U|$, then $X_{R,(H,U)}(i)\le |U|!n^\ell$.\label{fact:extension type trivial bound}
\end{enumerate}
\end{fact}

\COMMENT{For \ref{fact:extension type l kappa}, apply definition with $U'=U$ and rearrange. For (ii), note that definition trivially holds if right hand side is $\le 0$. Since we choose the minimal $\kappa$, we have $\kappa \le |V(H)\sm U'| \le |V(H)\sm U|=\ell$. \ref{fact:extension type empty} follows from the previous if $|H|=0$.}

Recall that in our process, we have $p_{\cC}(i)=i/\binom{n}{3}$, and since $i=\O(n^2)$, we have $p_{\cC}(i)=\O(1/n)$. Fact~\ref{fact:extension type trivial}\ref{fact:extension type l kappa} tells us that if the triples in $\cC(i)$ appeared independently at random, then we would have $\expn {X_{R,(H,U)}(i)}=\O(n^{\kappa(H,U)})$.

In order to track $X_{R,(H,U)}(i)$ during the process, the following observations are crucial. We will use (i) to establish a trend hypothesis and (ii) to establish a boundedness hypothesis for $X_{R,(H,U)}$.

\begin{fact}\label{fact:extension induction}
Let $(H,U)$ be a $\kappa$-balanced extension type and $e\in H$. Then the following hold:
\begin{enumerate}[label=\rm{(\roman*)}]
\item $(H-e,U)$ is a $(\kappa+1)$-balanced extension type.
\item $(H-H[U\cup e],U\cup e)$ is a $\kappa$-balanced extension type.\COMMENT{could more generally replace $U\cup e$ by any $U'\supseteq U$}
\end{enumerate}
\end{fact}

\proof
(i) is obvious. For (ii), let $H_1:=H-H[U\cup e]$. Note that $(H_1,U\cup e)$ clearly is an extension type. For any $U\cup e \In U' \In V(H_1)=V(H)$ we have $|H-H[U']|\ge |V(H)\sm U'|-\kappa$ since $(H,U)$ is $\kappa$-bounded. As $H_1-H_1[U']=H-H[U']$, we conclude that $|H_1-H_1[U']|\ge |V(H_1)\sm U'|-\kappa$.
\endproof

We will now make some observations as to how \Erd-configurations (and combinations thereof) can be viewed as balanced extension types.

\begin{fact}\label{fact:simple Erdos extension}
Let $\cS$ be an \Erd-configuration with $T'\in\cS$ and let $U\In V(\cS)$ with $|U|\ge 4$. Let $\cS^-:=\cS-\Set{T'}$. Then $|\cS^--\cS^-[U]|\ge |V(\cS)\sm U|$.
\end{fact}

\proof
Suppose not. Then $|\cS^-[U]|> |U|+|\cS^-|-|V(\cS)|=|U|-3$. But this means that $\cS^-[U]$ contains a forbidden configuration, a contradiction.
\endproof

The following result will be used to establish various boundedness conditions.

\begin{prop}\label{prop:Erdos extension analysis}
Let $\cS$ be an \Erd-configuration on $j$ points and $\cS'\In \cS$ with $V(\cS')=V(\cS)$ and $|\cS'|=c$. Let $U\In V(\cS)$ with $|U|\ge 4$ and $|(\cS-\cS')[U]|\ge a$. Define $\kappa:=j-3-c-a$. Then $(\cS'-\cS'[U],U)$ is $\max\Set{\kappa,0}$-balanced.
\end{prop}

\proof
Suppose not. Then there exists $U\In U'\In V(\cS)$ such that $|\cS'-\cS'[U']|< |V(\cS)\sm U'|-\max\Set{\kappa,0}$. Note that this implies $U'\neq V(\cS)$ and that $|\cS'|-|\cS'[U']|\le |V(\cS)|-|U'|-\kappa-1$. Thus $|\cS'[U']|\ge |U'|-2-a$.\COMMENT{$|\cS'[U']|\ge |U'|+(j-3-c-a)+1-j+c$} We conclude that $|\cS[U']|\ge |\cS'[U']|+|(\cS-\cS')[U]| \ge |U'|-2$. Since $|U'|\ge 4$, $\cS[U']$ is a forbidden configuration, a contradiction to $\cS$ being an \Erd-configuration.
\endproof

We will also need to bound the number of pairs $\cS,\cS'$ of \Erd-configurations appearing in some specified constellation. 

The following proposition yields a `global' edge count of two overlapping \Erd-configurations. After specifying some root set, it can be used to compare the number of edges with the number of free vertices.
\begin{prop}\label{prop:conf overlap}
Let $\cS_1,\cS_2$ be distinct \Erd-configurations.
\begin{enumerate}[label=\rm{(\roman*)}]
\item If $|V(\cS_1)\cap V(\cS_2)|\ge 4$, then $|\cS_1\cup \cS_2| \ge |V(\cS_1)\cup V(\cS_2)|-1$.
\item If $|V(\cS_1)\cap V(\cS_2)|= 3$, then $|\cS_1\cup \cS_2| \ge |V(\cS_1)\cup V(\cS_2)|-2$.
\end{enumerate}
\end{prop}

\proof
(i) View $\cS_1\cap \cS_2$ as a set of triples on $V(\cS_1)\cap V(\cS_2)$. Since $\cS_2$ is an \Erd-configuration, we must have $|\cS_1\cap \cS_2|\le |V(\cS_1)\cap V(\cS_2)|-3$. This implies that
\begin{align*}
|\cS_1\cup \cS_2|&=|\cS_1|+|\cS_2|-|\cS_1\cap \cS_2|=|V(\cS_1)|+|V(\cS_2)|-4-|\cS_1\cap \cS_2| \ge |V(\cS_1)\cup V(\cS_2)|-1.
\end{align*}

(ii) follows similarly by using $|\cS_1\cap \cS_2|\le 1$.
\endproof

\begin{prop} \label{prop:deryk counting 2}
Let $\cS_1,\cS_2$ be distinct \Erd-configurations with $T'\in \cS_1\cap \cS_2$. Define $H:=(\cS_1\cup \cS_2)-\Set{T'}$. Suppose $U'\In V(H)$ is such that $|U'\cap V(\cS_1)|\ge 4$ and $|(U'\cup V(\cS_1))\cap V(\cS_2)|\ge 4$. Then $|H-H[U']|\ge |V(H)\sm U'|$.
\end{prop}
\proof
Let $U_1':=U'\cap V(\cS_1)$ and $U_2':=(U'\cup V(\cS_1))\cap V(\cS_2)$. For $\ell\in[2]$, let $\cS_\ell^-:=\cS_\ell-\Set{T'}$. By Fact~\ref{fact:simple Erdos extension}, we have $|\cS_1^- -\cS_1^-[U_1']|\ge |V(\cS_1)\sm U_1'|$ and $|\cS_2^- -\cS_2^-[U_2']|\ge |V(\cS_2)\sm U_2'|$. 
Observe that $\cS_1^- -\cS_1^-[U_1'], \cS_2^- -\cS_2^-[U_2']  \In  H-H[U']$ and $(\cS_1^- -\cS_1^-[U_1'])\cap (\cS_2^- -\cS_2^-[U_2'])=\emptyset$.\COMMENT{Suppose $e\in \cS_1^-\cap \cS_2^-$. Then $e\In V(\cS_1)\cap V(\cS_2)$ and thus $e\In U_2'$, hence $e\notin \cS_2^- -\cS_2^-[U_2']$.}
We conclude that
\begin{align*}
|H-H[U']|&\ge |\cS_1^- -\cS_1^-[U_1']|+|\cS_2^- -\cS_2^-[U_2']|\\
         &\ge |V(\cS_1)\sm U_1'|+|V(\cS_2)\sm U_2'| =|V(H) \sm U'|,
\end{align*}
as required.
\endproof

The next two propositions will be used to establish trend and boundedness hypotheses for our key variables.

\begin{prop} \label{prop:double extension}
Let $T_1,T_2,T'$ be distinct triples and let $\cS_1,\cS_2$ be distinct \Erd-configurations where at least one is not the diamond. Suppose that $T'\in \cS_1\cap \cS_2$, $T_1\in \cS_1$, $T_2\in \cS_2$.
Define $H:=(\cS_1\cup \cS_2)-\Set{T_1,T_2,T'}$ and $U:=T_1\cup T_2$. Suppose that $H[U]$ is empty. Then the extension type $(H,U)$ is $0$-balanced. 
\end{prop}

\proof
Note that $|H|=|\cS_1\cup \cS_2|-3$ and $V(H)=V(\cS_1)\cup V(\cS_2)$.
By Proposition~\ref{prop:conf overlap}, we have $|\cS_1\cup \cS_2| \ge |V(\cS_1)\cup V(\cS_2)|-2$, implying $|H|\ge |V(H)|-5$. Thus, if $|U|\ge 5$, we have $|H|\ge |V(H)\sm U|$. Suppose $|U|=4$, i.e.~$|T_1\cap T_2|=2$. Since either $T_1$ and $T'$ or $T_2$ and $T'$ are edge-disjoint,\COMMENT{one could be a diamond} we have $T_1\cap T_2\not\In T'$, implying $|V(\cS_1)\cap V(\cS_2)|\ge 4$. Proposition~\ref{prop:conf overlap}(i) implies $|\cS_1\cup \cS_2| \ge |V(\cS_1)\cup V(\cS_2)|-1$ and thus $|H|\ge |V(H)|-4=|V(H)\sm U|$.

Now, let $U\In U' \In V(H)$ with $U'\neq U$. Without loss of generality we may assume that $|U'\cap V(\cS_1)|\ge 4$. Moreover, $|(U'\cup V(\cS_1))\cap V(\cS_2)|\ge |T_2\cup T'|\ge 4$. Thus, we can apply Proposition~\ref{prop:deryk counting 2} to conclude that $|H-H[U']|\ge |V(H)\sm U'|$.
\endproof

\begin{prop}\label{prop:butterfly count}
Let $\cS_1,\cS_2$ be distinct \Erd-configurations with distinct $T\in \cS_1$ and $T'\in \cS_1\cap \cS_2$ and such that $|V(\cS_1)\cap V(\cS_2)|\ge 4$. Define $H:=(\cS_1\cup \cS_2)-\Set{T,T'}$. Then $(H,T)$ is $0$-balanced.
\end{prop}
\proof
Clearly, $H[T]$ is empty.
Let $T\In U' \In V(H)$. Note that $|H|=|\cS_1\cup \cS_2|-2$ and $V(H)=V(\cS_1)\cup V(\cS_2)$. Since $|V(\cS_1)\cap V(\cS_2)|\ge 4$, we have $|\cS_1\cup \cS_2| \ge |V(\cS_1)\cup V(\cS_2)|-1$ by Proposition~\ref{prop:conf overlap}. We conclude that $|H|\ge |V(H)|-3=|V(H)\sm T|$. Thus, if $U'=T$, our balancedness requirement is satisfied.
Now, suppose $U'\neq T$. If $|U'\cap V(\cS_1)|\ge 4$ or $|U'\cap V(\cS_2)|\ge 4$, then Proposition~\ref{prop:deryk counting 2} (with the roles of $\cS_1$ and $\cS_2$ being swapped in the latter case) implies that $|H-H[U']|\ge |V(H)\sm U'|$. The only remaining case is if $U'\cap V(\cS_1)=T$ and $|U'\cap V(\cS_2)|\le 3$. Since $|\cS_2[U']|\le 1\le  |U'\sm T|$ and $|(\cS_1-T)[U']|=0$, we obtain $|H[U']|\le |U'\sm T|$, which together with $|H|\ge |V(H)\sm T|$ implies $|H-H[U']|\ge |V(H)\sm U'|$, as desired.
\endproof

\subsection{Stopping and freezing times}

In order to keep track of our key variables, we define the following error function:
\begin{align}
\eps(i) &:= \left(1+\frac{C}{n^2}\right)^i \eps_0. \label{error function}
\end{align}
Note that we have
\begin{align}
\Delta \eps(i)= C\eps(i)n^{-2} \mbox{ and } \eps(i)\le \eul^C \eps_0. \label{error function props}
\end{align}

To control the extension types, for all $\ell\in \Set{1,\dots,m}$, $\kappa \in\Set{0,\dots,\ell}$, define
\begin{align}
\eps_{\kappa,\ell}(i):= n^{\kappa+\frac{\ell}{m+\kappa}} (1+i/n^2). \label{extension error def}
\end{align}

Let $\tau_{cut}:=\lfloor (1-\gamma)n^2/6\rfloor $.
During the process (at least up to time $\tau_{cut}$), we aim to show that the following hold:
\begin{itemize}
\item for all $e\in E(i)$, we have
\begin{align}
X_e(i) &= f_{edge}(i) \pm \eps(i) n,      \label{edge regularity}
\end{align}
\item for all $T\in \cA(i)$, $j\in\Set{6,\dots,j_{max}}$ and $c\in \Set{0,\dots,j-4}$, we have
\begin{align}
X_{T,j,c}(i) &= f_{j,c}(i) \pm \eps(i)n^{j-3-c},    \label{erdos regularity}
\end{align}
\item for all $\ell\in \Set{1,\dots,m}$, $\kappa \in\Set{0,\dots,\ell}$, $(H,U)\in Ext(\kappa,\ell)$ and all sets $R\In V$ with $|R|=|U|$, we have
\begin{align}
X_{R,(H,U)}(i) &\le \eps_{\kappa,\ell}(i).   \label{extension bound}
\end{align}
\end{itemize}

Let $\tau_{violated}$ be the smallest $i$ such that at least one of these conditions is violated ($\tau_{violated}:=\infty$ if this never happens). Let
\begin{align}
\tau_{stop}:=\tau_{violated} \wedge \tau_{cut}. \label{def tau stop}
\end{align}
Clearly, the (random) times $\tau_{cut},\tau_{violated},\tau_{stop}$ are stopping times of the process (so for example whether $\tau_{violated}= i$ can be decided upon observing the process until time $i$). We now define additional `freezing times'. 

Define $\tau_{freeze,e}:= \tau_{stop}\wedge (\tau_e-1)$ and $\tau_{freeze,T}:= \tau_{stop}\wedge (\tau_T-1)$. (Recall that $\tau_e$ and $\tau_T$ were defined at the beginning of Section~\ref{sec:process}.) We note that the random times $\tau_{freeze,e}$ and $\tau_{freeze,T}$ are not stopping times of the process. 

For every $2$-set $e$, define 
\begin{align}
X_e^\pm(i) &:= \begin{cases}
\pm X_e(i) \mp f_{edge}(i) - \eps(i) n 	  &    \mbox{if } i\le \tau_{freeze,e}, \\   \label{margin def formal edge}     
X_e^\pm (i-1)       &    \mbox{if } i>\tau_{freeze,e}.
\end{cases}
\end{align}
Alternatively, we can write $$X_e^\pm(i) :=\pm X_e(i \wedge \tau_{freeze,e}) \mp f_{edge}(i\wedge \tau_{freeze,e}) - \eps(i\wedge \tau_{freeze,e}) n.$$

For every triple $T$, $j\in\Set{6,\dots,j_{max}}$ and $c\in\Set{0,\dots,j-4}$, define
\begin{align}
X_{T,j,c}^\pm(i) &:= \begin{cases}
\pm X_{T,j,c}(i) \mp f_{j,c}(i) - \eps(i) n^{j-3-c} 	  &   \mbox{if } i\le \tau_{freeze,T}, \\   \label{margin def formal triple}     
X_{T,j,c}^\pm (i-1)      &   \mbox{if }  i>\tau_{freeze,T}.
\end{cases}
\end{align}

For all $\ell\in \Set{1,\dots,m}$, $\kappa \in\Set{0,\dots,\ell}$, $(H,U)\in Ext(\kappa,\ell)$ and all sets $R\In V$ with $|R|=|U|$, define
\begin{align}
X^+_{R,(H,U)}(i) &:= \begin{cases}
X_{R,(H,U)}(i) - \eps_{\kappa,\ell}(i) 	  &  \mbox{if }  i\le \tau_{stop}, \\   \label{margin def formal extension}     
X^+_{R,(H,U)}(i-1)      &   \mbox{if }  i>\tau_{stop}.
\end{cases}
\end{align}

Recall from Section~\ref{sec:process} that $\tau_{max}:=\min\set{i}{\cA(i)=\emptyset}$.

\begin{fact}
Suppose that all the variables $X^{\pm}$ defined in \eqref{margin def formal edge}, \eqref{margin def formal triple}, \eqref{margin def formal extension} are non-positive for all $i$. Then $\tau_{max}\ge \tau_{cut}$.
\end{fact}

\proof
Clearly, if $E(\tau_{max})=\emptyset$, then $\tau_{max}=\lfloor \binom{n}{2}/3 \rfloor$. Thus, we can assume that there exists $e^\ast\in E(\tau_{max})$.

We first claim that $\tau_{violated}\ge \tau_{cut}$. Suppose for a contradiction that this is not the case. Say, for example, that condition~\eqref{edge regularity} for $e$ is violated at time $\tau_{violated}$. In particular, we have $\tau_{violated}\le \tau_e-1$. We conclude that $\tau_{violated}=\tau_{freeze,e}$. However, $X^+_e(\tau_{violated})>0$ or $X^-_e(\tau_{violated})>0$ since \eqref{edge regularity} is violated for $e$ at time $\tau_{violated}$, a contradiction to our assumption. The argument for the case when \eqref{erdos regularity} or~\eqref{extension bound} is violated is similar.\COMMENT{If the condition for $T,j,c$ is violated at time $\tau_{violated}$, we have $\tau_{violated}\le \tau_T-1$. We conclude that $\tau_{violated}=\tau_{freeze,T}$. Hence, $X^+_{T,j,c}(\tau_{violated})>0$ or $X^-_{T,j,c}(\tau_{violated})>0$, a contradiction to our assumption. If the condition for $R,(H,U)$ is violated at time $\tau_{violated}$, we have $\tau_{violated}=\tau_{stop}$. Hence, $X^+_{R,(H,U)}(\tau_{violated})>0$, a contradiction to our assumption.}

We deduce that $\tau_{stop}= \tau_{cut}$ and hence $\tau_{freeze,e^\ast}=\tau_{cut}$.\COMMENT{Since $e^\ast\in E(\tau_{max})$, we have $\tau_{e^\ast}=\infty>\tau_{stop}$. Thus $\tau_{freeze,e^\ast}=\tau_{stop}$.} From $X^-_{e^\ast}(\tau_{cut}) \le 0$, we infer that $|\cX_{e^\ast}(\tau_{cut})|=X_{e^\ast}(\tau_{cut}) \ge f_{edge}(\tau_{cut})-\eps(\tau_{cut}) n  >0$, where the last inequality follows from \eqref{def densities}, \eqref{edge trajectory}, \eqref{rho magnitudes} and~\eqref{error function props}. In particular, $\cA(\tau_{cut})\neq \emptyset$ and hence $\tau_{max}\ge \tau_{cut}$.
\endproof

The following lemma thus implies Theorem~\ref{thm:process stop}.

\begin{lemma}\label{lem:negative trend variables}
Whp, all the variables $X^{\pm}$ defined in \eqref{margin def formal edge}, \eqref{margin def formal triple}, \eqref{margin def formal extension} are non-positive for all $i$.
\end{lemma}

We first remark that all these variables are negative at the start. Since $\cA(0)=\binom{V}{3}$ and $\cC(0)=\emptyset$, we observe that for every $2$-set $e$, we have $X_e(0)=n-2=f_{edge}(0)$, and for every triple $T$, we have $X_{T,j,c}(0)=0=f_{j,c}(0)$ if $c>0$ and $X_{T,j,0}(0)=J_j=f_{j,0}(0)$. Moreover, for every extension type $(H,U)$ with $H$ not being empty, we have $X_{R,(H,U)}(0)=0$. Hence, by \eqref{margin def formal edge}--\eqref{margin def formal extension}, the following initial conditions hold for our variables $X^\pm$:
\begin{align}
X^\pm_e(0) &= -\eps_0 n;  \label{initial edge} \\
X^\pm_{T,j,c}(0) &= -\eps_0 n^{j-3-c};    \label{initial triple} \\
X^+_{R,(H,U)}(0) &= -n^{\kappa+\frac{\ell}{m+\kappa}}. \label{initial extension}
\end{align}

Our strategy to prove Lemma~\ref{lem:negative trend variables} is as follows: In the next subsection, we show that each such variable $X^\pm$ induces a supermartingale. We then establish some additional boundedness conditions which we need to finally apply Freedman's inequality to prove Lemma~\ref{lem:negative trend variables}.

\subsection{Counting double configurations}
For a triple $T$, we also define the variable $X_{T,double}(i)$ which counts the number of pairs $\cS,\cS'\in \bigcup_{j=4}^{j_{max}}\cX_{T,j,j-4}(i)$ such that $\cS\neq \cS'$ and $(\cS-\Set{T})\cap \cA(i)=(\cS'-\Set{T})\cap \cA(i)$.
Recall that $\cT_T(i)$ was defined in the beginning of Section~\ref{subsec:key vars}.

\begin{prop}\label{prop:danger count}
For all $i$ and $T\in \cA(i)$, $0\le  \sum_{j=4}^{j_{max}}X_{T,j,j-4}(i) - |\cT_T(i)| \le 2 X_{T,double}(i)$.
\end{prop}

\proof
For $T^\ast\in \cT_T(i)$, let $z_{T^\ast}$ denote the number of $\cS\in \bigcup_{j=4}^{j_{max}}\cX_{T,j,j-4}(i)$ with $(\cS-\Set{T})\cap \cA(i)=\Set{T^\ast}$. Thus, by Fact~\ref{fact:threats via extensions}, we have $\sum_{j=4}^{j_{max}}X_{T,j,j-4}(i)= \sum_{T^\ast\in \cT_T(i)} z_{T^\ast} \ge |\cT_T(i)|$, which establishes the first inequality. Crucially, we have $\sum_{T^\ast\in \cT_T(i)}\binom{z_{T^\ast}}{2}= X_{T,double}(i)$, implying that $\sum_{T^\ast\in \cT_T(i)\colon z_{T^\ast}>1}z_{T^\ast} \le 2 X_{T,double}(i)$. Thus,%
\COMMENT{For the inequality in the 1st line we use that $z_{T^*}\ge 1$ for each $T^*\in \cT_{T}(i)$ by Fact~\ref{fact:threats via extensions}.}
\begin{align*}
|\cT_T(i)| &\ge |\set{T^\ast\in \cT_T(i)}{z_{T^\ast}=1}|=\sum_{T^\ast\in \cT_T(i)}z_{T^\ast}-\sum_{T^\ast\in \cT_T(i) \colon z_{T^\ast}>1}z_{T^\ast}\\
								&\ge \sum_{j=4}^{j_{max}}X_{T,j,j-4}(i) -2 X_{T,double}(i).
\end{align*}
\endproof

The following is an immediate consequence of Proposition~\ref{prop:butterfly count}.

\begin{cor} \label{cor:butterfly count}
For all triples $T$ and all times $i< \tau_{stop}$, $X_{T,double}(i)=\O(n^{1-\frac{1}{m}})$.
\end{cor}

\proof
Suppose the pair $\cS,\cS'$ is counted by $X_{T,double}(i)$. Let $T'\in \cA(i)$ be the unique available triple in $\cS-\Set{T}$ and $\cS'-\Set{T}$. Let $H:=(\cS\cup \cS')-\Set{T,T'}$. Note that $H\In \cC(i)$. By Proposition~\ref{prop:butterfly count}, the extension type $(H,T)$ is $0$-balanced. Thus, $(H,T)\in Ext(0,\ell)$%
\COMMENT{Since  $\cS\neq \cS'$ we must have $H\neq \emptyset$.}
with $\ell:=|V(H)\sm T|\le m-1$ (with room to spare). 

Hence, we have $$X_{T,double}(i)= \O\left( \sum_{\ell=0}^{m-1}\sum_{(H',T)\in Ext(0,\ell)} X_{T,(H',T)}(i)\right) .$$
By~\eqref{extension bound} and \eqref{extension error def}, we conclude that $X_{T,double}(i) =\O(n^{\frac{m-1}{m}})$.
\endproof

For distinct triples $T_1,T_2$, we let $\cX_{T_1,T_2}(i)$ be the set of all pairs $\cS_1\neq \cS_2$, not both diamonds, such that for each $\ell\in [2]$, $\cS_\ell \in \cX_{T_\ell,j_\ell,j_\ell-4}(i)$ with $4\le j_\ell\le j_{max}$, and such that $(\cS_1-\Set{T_1})\cap \cA(i)=(\cS_2-\Set{T_2})\cap \cA(i)$. We let $X_{T_1,T_2}(i):=|\cX_{T_1,T_2}(i)|$.

Recall that $th_{\cS,T}(i)$ was defined after \eqref{threat approximation edge}.

\begin{prop} \label{prop:Erdos threats}
For all $i$, all $T\in \cA(i)$, all $j\in \Set{6,\dots,j_{max}}$, all $c\in \Set{0,\dots,j-4}$ and all $\cS\in \cX_{T,j,c}(i)$, $$\left|th_{\cS,T}(i)-\sum_{T'\in (\cS-\Set{T})\cap \cA(i)}|\cT_{T'}(i)|\right|\le \O(1)+\sum_{T'\neq T''\in \cS\cap \cA(i)} X_{T',T''}(i).$$
\end{prop}

\proof
We have
\begin{align*}
th_{\cS,T}(i) &= \left|\left(\bigcup_{T'\in (\cS-\Set{T})\cap \cA(i)} (\Set{T'}\cup \cT_{T'}(i)) \right) \sm (\Set{T}\cup \cT_T(i))\right|,
\end{align*}
from which we immediately have $th_{\cS,T}(i)\le j+\sum_{T'\in (\cS-\Set{T})\cap \cA(i)}|\cT_{T'}(i)|$. Moreover, we have
\begin{align*}
th_{\cS,T}(i) &\ge \sum_{T'\in (\cS-\Set{T})\cap \cA(i)}|\cT_{T'}(i)| - \sum_{T'\neq T''\in \cS\cap \cA(i)} |\cT_{T'}(i)\cap \cT_{T''}(i)| -1.
\end{align*}
We claim that $|\cT_{T'}(i)\cap \cT_{T''}(i)|\le X_{T',T''}(i) + \O(1)$, which completes the proof. Let $T^\ast \in \cT_{T'}(i)\cap \cT_{T''}(i)$. By Fact~\ref{fact:threats via extensions}, there are $\cS' \in \cX_{T',j',j'-4}(i)$ and $\cS'' \in \cX_{T'',j'',j''-4}(i)$ with $4\le j',j''\le j_{max}$ and such that $(\cS'-\Set{T'})\cap \cA(i)=\Set{T^\ast}=(\cS''-\Set{T''})\cap \cA(i)$. Clearly, we have $\cS'\neq \cS''$. Thus, unless $j'=j''=4$, this pair $\cS',\cS''$ is counted by $X_{T',T''}(i)$. Finally, if both $\cS',\cS''$ are diamonds, then since $T'$ and $T''$ are edge-disjoint,\COMMENT{since $j\ge 6$} we must have $T^\ast\In T'\cup T''$, for which there are only $\O(1)$ possibilities.
\endproof

The following fact will be useful in Section~\ref{subsec:boundedness} to bound the negative change of $X_{T,j,j-4}(i)$.

\begin{fact}\label{fact:negative danger bound}
Let $j\in \Set{6,\dots,j_{max}}$ and $T,T^\ast\in \cA(i)$ be distinct. Then
$$|\set{\cS\in \cX_{T,j,j-4}(i)}{T^\ast\mbox{ threatens }\cS,T}|\le X_{T,T^\ast}(i).$$ 
\end{fact}

\proof
Let $\cS\in \cX_{T,j,j-4}(i)$ and assume that $T^\ast$ threatens $\cS,T$. Let $\Set{T'}=(\cS-\Set{T})\cap \cA(i)$. We cannot have $T^\ast=T'$ as this would mean $T^\ast\leftrightarrow T$. Hence, by Fact~\ref{fact:threats via extensions}, there is $\cS' \in \cX_{T^\ast,j',j'-4}(i)$ with $4\le j'\le j_{max}$ such that $(\cS'-\Set{T^\ast})\cap \cA(i)=\Set{T'}$. Since $T\neq T^\ast$, the pair $\cS,\cS'$ is counted by $X_{T,T^\ast}(i)$.
\endproof

\begin{cor} \label{cor:double extension}
For $i< \tau_{stop}$ and distinct $T_1,T_2\in \cA(i)$, we have $X_{T_1,T_2}(i)=\O(n^{1-\frac{1}{m}})$.
\end{cor}

\proof
Let $U:=T_1\cup T_2$.
Suppose the pair $\cS_1,\cS_2$ is counted by $X_{T_1,T_2}(i)$. Let $T'$ be the unique available triple in $\cS_1-\Set{T_1}$ and $\cS_2-\Set{T_2}$. Let $H:=(\cS_1\cup \cS_2)-\Set{T_1,T_2,T'}$. Note that $H\In \cC(i)$. In particular, since $T_1,T_2$ are still available, we have that $H[U]$ is empty. By Proposition~\ref{prop:double extension}, the extension type $(H,U)$ is $0$-balanced. Thus, $(H,U)\in Ext(0,\ell)$%
\COMMENT{Since at most one of $\cS_1,\cS_2$ is a diamond, we must have $H\neq\emptyset$.}
with $\ell:=|V(H)\sm U|\le m-1$ (with room to spare). 

Hence, we have $$X_{T_1,T_2}(i)= \O\left(\sum_{\ell=0}^{m-1}\sum_{(H',U)\in Ext(0,\ell)} X_{U,(H',U)}(i)\right) .$$
By~\eqref{extension bound} and \eqref{extension error def}, we conclude that $X_{T_1,T_2}(i) =\O(n^{\frac{m-1}{m}})$.
\endproof

We will also need the following consequence of Proposition~\ref{prop:Erdos extension analysis}.

\begin{cor}\label{cor:Erdos extension applicable}
Let $R\In V$ with $|R|\ge 4$, $j\in \Set{6,\dots,j_{max}}$ and $c\in \Set{0,\dots,j-4}$. For $i<\tau_{stop}$, the number of \Erd-configurations $\cS$ with $R\In V(\cS)$, $|V(\cS)|=j$, $|\cS\cap \cC(i)|=c$ and $|(\cS-\cC(i))[R]|\ge a$, is $\O(n^{\max\Set{j-3-c-a,0} +\frac{1}{2}})$.
\end{cor}

\proof
Let $\cS$ be an \Erd-configuration with $R\In V(\cS)$, $|V(\cS)|=j$, $|\cS\cap \cC(i)|=c$ and $|(\cS-\cC(i))[R]|\ge a$. Let $\cS':=\cS\cap \cC(i)$ and $\kappa:=\max\Set{j-3-c-a,0}$. By Proposition~\ref{prop:Erdos extension analysis}\COMMENT{with $U=R$}, we know that $(\cS'-\cS'[R],R)$ is $\kappa$-balanced. Thus, unless $\cS'-\cS'[R]$ is empty, we have $(\cS'-\cS'[R],R)\in Ext(\kappa',\ell)$ for some $0\le \kappa'\le \kappa$ and $\ell:=j-|R|\le j_{max}= m/2$. Therefore, the number of such $\cS$ is $\O(n^{\kappa+\frac{j_{max}}{m}})$ by~\eqref{extension bound} and \eqref{extension error def}. In case $\cS'-\cS'[R]$ is empty, we have that the number of free vertices is $j-|R|\le \kappa$ by Fact~\ref{fact:extension type trivial}\ref{fact:extension type empty}, and thus obtain the trivial upper bound $\O(n^{\kappa})$.
\endproof

\subsection{Trend hypotheses}

Our goal is now to show that the variables $X^\pm$ form supermartingales.
For $i\ge 0$, define the random variable $$L(i):=(T^\ast(0),T^\ast(1),\dots,T^\ast((i-1)\wedge (\tau_{max}-1)))$$ which lists the sequence of chosen triples until time~$i$. Thus, $L(i)$ contains all the information about the process until time $i$.
Let $\cL(i)$ denote the set of all possible outcomes of $L(i)$. Moreover, let $\cL^\ast(i)$ denote the set of all $\tilde{L}\in \cL(i)$ for which $\tau_{stop}>i$.

Now, let $X$ be one of our variables $X^\pm$ defined in \eqref{margin def formal edge}--\eqref{margin def formal extension}. We will show that $(X(0),X(1),\dots)$ is a supermartingale with respect to (the filtration induced by) $(L(0),L(1),\dots)$.
Thus, we need to show that for all $i \ge 0$ and all $\tilde{L}\in \cL(i)$, we have $$\expn{\Delta X(i)\mid L(i)=\tilde{L}}\le 0.$$
Recall that $X$ comes with a (random) `freezing' time $\tau$ (e.g.~$\tau_{freeze,e}$, $\tau_{freeze,T}$, $\tau_{stop}$), of which we know that 
\begin{align}
	\tau \le \tau_{stop} \text{ and $\Delta X(i)=0$ for all $i\ge \tau$}. \label{abstract freezing}
\end{align}

Consider $i\ge 0$ and $\tilde{L}\in \cL(i)$. We may transition to the probability space $\mathbb{P}_{\tilde{L}}$ obtained by conditioning on the event $L(i)=\tilde{L}$. Thus, we need to show that $\cexpn{\tilde{L}}{\Delta X(i)}\le 0$.
If $\cprob{\tilde{L}}{i<\tau} =0$, then trivially $\cexpn{\tilde{L}}{\Delta X(i)}= 0$ by~\eqref{abstract freezing}. 
Note that if $\tilde{L}\in \cL(i)\sm \cL^\ast(i)$, then we have $\cprob{\tilde{L}}{i<\tau} =0$ by~\eqref{abstract freezing}.
If $\tilde{L}\in \cL^\ast(i)$ and $\cprob{\tilde{L}}{i<\tau} >0$, then by the law of total expectation, we obtain
\begin{eqnarray*}
\cexpn{\tilde{L}}{\Delta X(i)} &=&  \cexpn{\tilde{L}}{\Delta X(i) \mid i<\tau } \cprob{\tilde{L}}{i<\tau} + \cexpn{\tilde{L}}{\Delta X(i) \mid i\ge\tau } \cprob{\tilde{L}}{i\ge\tau},
\end{eqnarray*}
where the second summand trivially vanishes, again by~\eqref{abstract freezing}.

To summarise, in order to show that $(X(0),X(1),\dots)$ is a supermartingale, it suffices to show that
\begin{align}
	\cexpn{\tilde{L}}{\Delta X(i) \mid i<\tau } \le 0   \label{supermartingale sufficient}
\end{align}
for all $i\ge 0$ and all $\tilde{L}\in \cL^\ast(i)$ with $\cprob{\tilde{L}}{i<\tau} >0$.

Similarly, in order to show that $\expn{|\Delta X(i)| \mid L(i)} \le K$, it suffices to show that 
\begin{align}
	\cexpn{\tilde{L}}{|\Delta X(i)| \mid i<\tau } \le K   \label{expected abs sufficient}
\end{align}
for all $i\ge 0$ and all $\tilde{L}\in \cL^\ast(i)$ with $\cprob{\tilde{L}}{i<\tau} >0$.

\begin{lemma}\label{lem:supermartingales}
The following hold:
\begin{enumerate}[label=\rm{(\roman*)}]
\item For every $2$-set $e$, $(X_e^+(0),X_e^+(1),\dots)$ and $(X_e^-(0),X_e^-(1),\dots)$ are supermartingales with respect to $(L(0),L(1),\dots)$. Moreover, $$\expn{|\Delta X^\pm_e(i)| \mid  L(i)}=\O_{\gamma}(n^{-1})$$ for all $i$.
\item For every triple $T$, all $j\in \Set{6,\dots,j_{max}}$ and $c\in \Set{0,\dots,j-4}$, $(X_{T,j,c}^+(0),X_{T,j,c}^+(1),\dots)$ and $(X_{T,j,c}^-(0),X_{T,j,c}^-(1),\dots)$ are supermartingales with respect to $(L(0),L(1),\dots)$. Moreover, $$\expn{|\Delta X^\pm_{T,j,c}(i)| \mid  L(i)}=\O_{\gamma}(n^{j-5-c})$$ for all $i$.
\item For all $\ell\in \Set{1,\dots,m}$, $\kappa \in\Set{0,\dots,\ell}$, $(H,U)\in Ext(\kappa,\ell)$ and all sets $R\In V$ with $|R|=|U|$, $(X_{R,(H,U)}^+(0),X_{R,(H,U)}^+(1),\dots)$ is a supermartingale with respect to $(L(0),L(1),\dots)$. Moreover, $$\expn{|\Delta X^+_{R,(H,U)}(i)| \mid  L(i)}\le 2n^{\kappa-2+\frac{\ell}{m+\kappa}}$$ for all $i$.
\end{enumerate}
\end{lemma}

We will prove this lemma at the end of this subsection.
To continue, we need to gain control over $|\cA(i)|$, $th_{T,e}(i)$ and $th_{\cS,T}(i)$.
Recall from~\eqref{edge trajectory}, \eqref{available trajectory} and~\eqref{def sum F} that 
\begin{align*}
A(i)&=\eul^{-\rho(i)}p(i)^3 \binom{n}{3}=\frac{1}{3}p(i)\binom{n}{2}f_{edge}(i), \\
F(i) &=\sum_{j=6}^{j_{max}} f_{j,j-4}(i).
\end{align*}

Note that if $i\le \tau_{cut}$, then 
\begin{align}
p(i)=\Omega_\gamma(1), \; A(i) &= \Omega_\gamma (n^3), \label{constant density cutoff}
\end{align}
where we use \eqref{rho magnitudes} to deduce the latter from the first.

\begin{lemma} \label{lem:properties before stopping} 
Let $i<\tau_{stop}$. Then the following hold:
\begin{enumerate}[label=\rm{(\roman*)}]
\item $|\cA(i)|= A(i) \pm \eps(i) n^3$.\label{properties before stopping number available}
\item For all $T\in \cA(i)$, $|\cT_T(i)|= \O(n)$. \label{properties before stopping number threats}
\item For all $e\in E(i)$ and $T\in \cX_e(i)$, $$th_{T,e}(i)= 2 f_{edge}(i) + F(i) + \O(\eps(i))n.$$\label{properties before stopping number edge threats}
\item For all $T\in \cA(i)$, all $j\in \Set{6,\dots,j_{max}}$, $c\in \Set{0,\dots,j-4}$ and $\cS\in \cX_{T,j,c}(i)$, $$th_{\cS,T}(i) = (j-3-c)(3 f_{edge}(i) + F(i)) + \O(\eps(i))n.$$\label{properties before stopping number triple threats}
\end{enumerate}
\end{lemma}

\proof
By~\eqref{def tau stop}, $i<\tau_{violated}$, so we can make use of~\eqref{edge regularity},~\eqref{erdos regularity} and~\eqref{extension bound}. Using Fact~\ref{fact:available edge triple relation}, we can deduce that
\begin{align*}
|\cA(i)|= \frac{1}{3}\sum_{e\in E(i)} X_e(i) \overset{\eqref{def densities},\eqref{edge regularity}}{=} \frac{1}{3}p(i)\binom{n}{2} (f_{edge}(i)\pm \eps(i) n) = A(i) \pm \eps(i) n^3,
\end{align*}
i.e.~\ref{properties before stopping number available} holds.

Moreover, from Fact~\ref{fact:smallest Erdos},~\eqref{edge regularity},~\eqref{erdos regularity} and~\eqref{def sum F}, we obtain that for all $T\in \cA(i)$ we have
\begin{align*}
\sum_{j=4}^{j_{max}}X_{T,j,j-4}(i)  = 3 f_{edge}(i) + F(i) + \O(\eps(i))n.
\end{align*}
From Corollary~\ref{cor:butterfly count}, we infer that $X_{T,double}(i)\le \eps(i) n$ for all triples $T$. With Proposition~\ref{prop:danger count} and the above, we have
\begin{align}
|\cT_T(i)| &= 3 f_{edge}(i) + F(i) + \O(\eps(i))n.\label{triple threat count}
\end{align}
 for all $T\in \cA(i)$, so~\ref{properties before stopping number threats} holds by~\eqref{edge and danger magnitude}.

By~\eqref{triple threat count}, Proposition~\ref{prop:edge triple dangers} and~\eqref{edge regularity}, for all $e\in E(i)$ and $T\in \cX_e(i)$ we have $$th_{T,e}(i)= 2 f_{edge}(i) + F(i) + \O(\eps(i))n,$$ i.e.~\ref{properties before stopping number edge threats} holds.

Finally, let $T\in \cA(i)$, $j\in \Set{6,\dots,j_{max}}$, $c\in \Set{0,\dots,j-4}$ and $\cS\in \cX_{T,j,c}(i)$. From Proposition~\ref{prop:Erdos threats} and Corollary~\ref{cor:double extension}, we deduce that
\begin{align*}
th_{\cS,T}(i) = \sum_{T'\in (\cS-\Set{T})\cap \cA(i)}|\cT_{T'}(i)| + \O(n^{1-\frac{1}{m}}) \overset{\eqref{triple threat count}}{=} (j-3-c)(3 f_{edge}(i) + F(i)) + \O(\eps(i))n.
\end{align*} Thus \ref{properties before stopping number triple threats} holds too.
\endproof

In order to compare the expectation of $\Delta X$ with its trajectory, we collect some important properties of the relevant trajectories in the following lemma.

\begin{lemma} \label{lem:derivatives}
The following hold for $0\le i\le \tau_{cut}$:
\begin{enumerate}[label=\rm{(\roman*)}]
\item $f'_{edge}(i) = -\frac{(2 f_{edge}(i) + F(i))f_{edge}(i)}{A(i)}$. \label{lem:derivatives edge}
\item For all $j\in\Set{6,\dots,j_{max}}$ and $c\in \Set{1,\dots,j-4}$, $$f'_{j,c}(i)  =  \frac{-(j-3-c)(3 f_{edge}(i) + F(i))f_{j,c}(i)  +  (j-2-c)f_{j,c-1}(i)}{A(i)}.$$ \label{lem:derivatives erdos}
\item For all $j\in\Set{6,\dots,j_{max}}$, $f'_{j,0}(i) = \frac{-(j-3)(3 f_{edge}(i) + F(i))f_{j,0}(i)}{A(i)}$. \label{lem:derivatives erdos zero}
\item $f'_{edge}(i)=\O_\gamma(n^{-1})$ and $f''_{edge}(i) = \O_{\gamma}(n^{-3})$. \label{lem:derivatives edge magnitude}
\item For all $j\in\Set{6,\dots,j_{max}}$ and $c\in \Set{0,\dots,j-4}$, $f'_{j,c}(i)=\O_\gamma(n^{j-5-c})$ and $f''_{j,c}(i) = \O_{\gamma}(n^{j-7-c})$. \label{lem:derivatives erdos magnitude}
\item $\Delta f_{edge}(i) = f'_{edge}(i) + \O_\gamma(n^{-3})= \O_\gamma(n^{-1})$.\label{edge trajectory taylor}
\item For all $j\in\Set{6,\dots,j_{max}}$ and $c\in \Set{0,\dots,j-4}$, $\Delta f_{j,c}(i)  = f'_{j,c}(i)+ \O_\gamma(n^{j-7-c}) = \O_\gamma(n^{j-5-c})$.\label{erdos trajectory taylor}
\end{enumerate}
\end{lemma}

\proof
First, we observe the following key identities:
\begin{eqnarray}
\frac{F(i)}{A(i)} &=&  \frac{\sum_{j=6}^{j_{max}} f_{j,j-4}(i)}{\eul^{-\rho(i)}p(i)^3 \binom{n}{3}}  \overset{\eqref{triple trajectory}}{=}  \sum_{j=6}^{j_{max}}\frac{(j-3)J_j}{\binom{n}{3}^{j-3}}i^{j-4} \overset{\eqref{def rho}}{=} \rho'(i),  \label{big F over A} \\
\frac{f_{edge}(i)}{A(i)} &\overset{\eqref{available trajectory}}{=}& \frac{6}{p(i)n(n-1)}  \overset{\eqref{edge density}}{=}  -\frac{p'(i)}{p(i)}, \label{small f over A}
\end{eqnarray}
and for all $j\in\Set{6,\dots,j_{max}}$ and $c\in \Set{1,\dots,j-4}$,
\begin{align}
\frac{f_{j,c-1}(i)}{f_{j,c}(i)} &\overset{\eqref{triple trajectory}}{=} \frac{\binom{j-3}{c-1}}{\binom{j-3}{c}}\eul^{-\rho(i)}p(i)^3 i^{-1}\binom{n}{3} = \frac{c}{(j-2-c)i}A(i). \label{small f c}
\end{align}

Using the chain rule, we can now easily check that
\begin{eqnarray*}
f'_{edge}(i) &=& -\rho'(i)f_{edge}(i) + 2\frac{p'(i)}{p(i)}f_{edge}(i) \overset{\eqref{big F over A},\eqref{small f over A}}{=} -\frac{(2 f_{edge}(i) + F(i))f_{edge}(i)}{A(i)},  \\
f'_{j,c}(i) &=& -(j-3-c) \rho'(i) f_{j,c}(i) + 3(j-3-c)\frac{p'(i)}{p(i)} f_{j,c}(i) + \frac{c}{i}f_{j,c}(i) \\
						&\overset{\eqref{big F over A},\eqref{small f over A}}{=}& \frac{-(j-3-c)(3 f_{edge}(i) + F(i))f_{j,c}(i)}{A(i)} +  \frac{c}{i}f_{j,c}(i),
\end{eqnarray*}
where the last summand vanishes if $c=0$ and can be replaced with $\frac{(j-2-c)f_{j,c-1}(i)}{A(i)}$ otherwise by~\eqref{small f c}. Hence, \ref{lem:derivatives edge}, \ref{lem:derivatives erdos} and \ref{lem:derivatives erdos zero} hold.

We continue with computing the second derivatives. Note that $\left(\frac{p'(i)}{p(i)}\right)'= - \frac{p'(i)^2}{p(i)^2}$. Therefore,
\begin{align*}
f''_{edge}(i) &= -\rho''(i)f_{edge}(i) - \rho'(i)f'_{edge}(i) - 2\frac{p'(i)^2}{p(i)^2}f_{edge}(i) + 2\frac{p'(i)}{p(i)}f'_{edge}(i), \\
f''_{j,c}(i) &= -(j-3-c)\left(f'_{j,c}(i) (-3 \frac{p'(i)}{p(i)} + \rho'(i)) + f_{j,c}(i) (3 \frac{p'(i)^2}{p(i)^2} + \rho''(i))  \right)       \\
              &+ (j-2-c) \frac{f'_{j,c-1}(i)A(i)-f_{j,c-1}A'(i)}{A(i)^2},
\end{align*}
where the last summand is not present if $c=0$. We clearly have $p'(i)=\O(n^{-2})$. Moreover, for the specified range of $i$, we have $\rho'(i)=\O(n^{-2})$ and $\rho''(i)=\O(n^{-4})$ by~\eqref{rho magnitudes} and, crucially, $p(i)=\Omega_\gamma(1)$ and $A(i)=\Omega_\gamma(n^3)$ by \eqref{constant density cutoff}. This also implies that $A'(i)=-\rho'(i)A(i)+3\frac{p'(i)}{p(i)}A(i)=\O_{\gamma}(n)$.

Together with~\eqref{triple trajectory magnitude}, \eqref{edge and danger magnitude},
we can infer that $f'_{edge}(i)=\O_\gamma(n^{-1})$ and $f'_{j,c}(i)=\O_\gamma(n^{j-5-c})$ and can conclude that $f''_{edge}(i)=\O_\gamma(n^{-3})$ and $f''_{j,c}(i)=\O_\gamma(n^{j-7-c})$. Thus, \ref{lem:derivatives edge magnitude} and \ref{lem:derivatives erdos magnitude} hold as well.

Finally, \ref{edge trajectory taylor} and \ref{erdos trajectory taylor} follow from the previous and \eqref{Taylor}.
\endproof

We are now in a position to show that the variables $X^\pm$ indeed form supermartingales.

\lateproof{Lemma~\ref{lem:supermartingales}}
Consider any $i\ge 0$ and any $\tilde{L}\in \cL^\ast(i)$. We consider the probability space $\mathbb{P}_{\tilde{L}}$. 
The values of all (random) variables at time $i$ are now determined by~$\tilde{L}$. Recall that by definition of $\cL^\ast(i)$, we have $i<\tau_{stop}$.
Hence, by~\eqref{constant density cutoff}, we have that $p(i)=\Omega_\gamma(1)$ and $A(i) = \Omega_\gamma (n^3)$.

\begin{step}
The expected change of $X_e$
\end{step}

Consider a $2$-set $e$. By the observation at~\eqref{supermartingale sufficient}, we may assume that $\cprob{\tilde{L}}{i<\tau_{freeze,e}}>0$. In particular, we have $e\in E(i)$. 

For every $T\in \cX_e(i)$, the probability that $T\notin \cX_e(i+1)$, conditioned on the event $e\in E(i+1)$, is $\frac{th_{T,e}(i)}{|\cA(i)\sm \cX_e(i)|}$ (cf.~Section~\ref{subsec:heuristics}). Thus,
\begin{align*}
\cexpn{\tilde{L}}{\Delta X_e(i) \mid i<\tau_{freeze,e}} = \cexpn{\tilde{L}}{\Delta X_e(i) \mid e\in E(i+1)}= - \sum_{T \in \cX_e(i)}\frac{th_{T,e}(i)}{|\cA(i)\sm \cX_e(i)|}.														
\end{align*}
By \eqref{edge regularity}, we have $|\cX_e(i)| = f_{edge}(i) \pm \eps(i) n$. By Lemma~\ref{lem:properties before stopping}\ref{properties before stopping number available}, we have $|\cA(i)|=A(i) \pm \eps(i) n^3$ and thus $|\cA(i)\sm \cX_e(i)|=A(i) \pm 2\eps(i) n^3$. Moreover, by Lemma~\ref{lem:properties before stopping}\ref{properties before stopping number edge threats}, we have $th_{T,e}(i)=2 f_{edge}(i) + F(i) + \O(\eps(i))n$ for all $T\in \cX_e(i)$. We conclude that 
\begin{eqnarray}\label{expected edge trend}
\cexpn{\tilde{L}}{\Delta X_e(i) \mid i<\tau_{freeze,e}} &=& - \left(f_{edge}(i) \pm \eps(i) n \right)\frac{2 f_{edge}(i) + F(i) + \O(\eps(i))n}{A(i) \pm 2\eps(i) n^3}  \nonumber  \\
																&\overset{\eqref{O fractions},\eqref{edge and danger magnitude}}{=}& -\frac{(2 f_{edge}(i) + F(i))f_{edge}(i)}{A(i)} + \O_{\gamma}(\eps(i)n^{-1})   \\
																&=& f'_{edge}(i) + \O_{\gamma}(\eps(i)n^{-1})  = \Delta f_{edge}(i) + \O_{\gamma}(\eps(i)n^{-1}), \nonumber
\end{eqnarray}
where the last two equalities are implied by Lemma~\ref{lem:derivatives}\ref{lem:derivatives edge} and~\ref{edge trajectory taylor}.

We conclude that
\begin{eqnarray*}
\cexpn{\tilde{L}}{\Delta X^{\pm}_e(i) \mid i<\tau_{freeze,e}} &\overset{\eqref{margin def formal edge}}{=}& \pm \cexpn{\tilde{L}}{\Delta X_e(i) \mid i<\tau_{freeze,e}} \mp \Delta f_{edge}(i) \\ 
    &  & -\Delta \eps(i) n  \\
            &\overset{\eqref{expected edge trend},\eqref{error function props}}{=}& \O_{\gamma}(\eps(i))n^{-1} -  C \eps(i)n^{-1} \overset{\eqref{hierarchy}}{\le} 0.												
\end{eqnarray*}
With the observation at \eqref{supermartingale sufficient}, this completes the proof that $(X_e^\pm(0),X_e^\pm(1),\dots)$ is a supermartingale with respect to $(L(0),L(1),\dots)$.

Moreover, since $\Delta X_e(i)\le 0$, we can also deduce that 
\begin{align}
\cexpn{\tilde{L}}{|\Delta X_e(i)| \mid i<\tau_{freeze,e}}=|\Delta f_{edge}(i)|+\O_{\gamma}(\eps(i)n^{-1})=\O_\gamma (n^{-1}) \label{abs bound edge}
\end{align}
by Lemma~\ref{lem:derivatives}\ref{lem:derivatives edge magnitude}.

\begin{step}
The expected change of $X_{T,j,c}$
\end{step}

Now, consider a triple $T$, $j\in \Set{6,\dots,j_{max}}$ and $c\in \Set{0,\dots,j-4}$.
By the observation at~\eqref{supermartingale sufficient}, we may assume that $\cprob{\tilde{L}}{i<\tau_{freeze,T}}>0$. In particular, we have $T\in \cA(i)$. 
We split the expected change of $X_{T,j,c}(i)$ into an expected loss and an expected gain, i.e.~$$\cexpn{\tilde{L}}{\Delta X_{T,j,c}(i) \mid i<\tau_{freeze,T}}= \cexpn{\tilde{L}}{\Delta X_{T,j,c}(i) \mid T\in \cA(i+1)}= -E^{loss} +E^{gain},$$ 
where $E^{loss}$ is the conditional expected size of $\cX_{T,j,c}(i)\sm \cX_{T,j,c}(i+1)$ and $E^{gain}$ is the conditional expected size of $\cX_{T,j,c}(i+1)\sm \cX_{T,j,c}(i)$.
As we condition on the event that $i<\tau_{freeze,T}$, the chosen triple $T^\ast(i)$ is chosen uniformly from the available triples which do not render $T$ unavailable, i.e.~$T^\ast(i)\notin \cT_T(i)\cup \Set{T}$.

We first consider $E^{loss}$.
For every $\cS \in \cX_{T,j,c}(i)$, we have that the (conditional) probability that $\cS \notin \cX_{T,j,c}(i+1)$, is $\frac{th_{\cS,T}(i)}{|\cA(i)\sm (\cT_T(i)\cup \Set{T})|}$ by definition of $th_{\cS,T}(i)$. By~\eqref{erdos regularity}, we have $|\cX_{T,j,c}(i)|=f_{j,c}(i) \pm \eps(i)n^{j-3-c}$. Thus, using Lemma~\ref{lem:properties before stopping}\ref{properties before stopping number available},\ref{properties before stopping number threats} and~\ref{properties before stopping number triple threats}, we conclude that
\begin{eqnarray*}
E^{loss} &=& \sum_{\cS \in \cX_{T,j,c}(i)} \frac{th_{\cS,T}(i)}{|\cA(i)\sm (\cT_T(i)\cup \Set{T})|} \\
         &=& \left( f_{j,c}(i) \pm \eps(i)n^{j-3-c} \right) \frac{(j-3-c)(3 f_{edge}(i) + F(i)) + \O(\eps(i))n}{A(i) \pm 2\eps(i) n^3} \\
         &\overset{\eqref{O fractions},\eqref{triple trajectory magnitude},\eqref{edge and danger magnitude}}{=}& \frac{(j-3-c)(3 f_{edge}(i) + F(i))}{A(i)}f_{j,c}(i) + \O_{\gamma}(\eps(i))n^{j-5-c}.
\end{eqnarray*}

We now consider $E^{gain}$.
Observe that if $\cS\in \cX_{T,j,c}(i+1)\sm \cX_{T,j,c}(i)$, then we must have $\cS\in \cX_{T,j,c-1}(i)$ and $T^\ast(i)\in \cS-\Set{T}$ (and in particular, $c>0$). Hence, if $c=0$, then $E^{gain}=0$. Assume now that $c>0$. For every $\cS\in \cX_{T,j,c-1}(i)$, we have $|(\cS-\Set{T})\cap \cA(i)|=j-3-(c-1)$ by definition of $\cX_{T,j,c-1}(i)$ (cf.~\eqref{def erdos}). Thus, there are $j-2-c$ available triples in $\cS-\Set{T}$ that, if chosen, could potentially imply $\cS\in \cX_{T,j,c}(i+1)$. However, some of these available triples might threaten $T$ or another available triple in~$\cS$. Thus, 
the (conditional) probability that $\cS \in \cX_{T,j,c}(i+1)$ is $$\frac{|((\cS-\Set{T})\cap \cA(i))\sm \bigcup_{T'\in \cS\cap \cA(i)}\cT_{T'}(i)|}{|\cA(i)\sm (\cT_T(i)\cup \Set{T})|}.$$
We claim that for most $\cS\in \cX_{T,j,c-1}(i)$, we have $(\cS-\Set{T})\cap \cA(i)\cap \bigcup_{T'\in \cS\cap \cA(i)}\cT_{T'}(i)=\emptyset$. Indeed, observe that the number of $\cS\in \cX_{T,j,c-1}(i)$ with $(\cS-\Set{T})\cap \cA(i)\cap \bigcup_{T'\in \cS\cap \cA(i)}\cT_{T'}(i)\neq \emptyset$ is at most the number of pairs $\cS,\cS'$ with $\cS\in \cX_{T,j,c-1}(i)$ and $\cS'\in \cX_{T',j',j'-4}(i)$, where $4\le j'\le j_{max}$, satisfying $\cS\cap \cS'\cap \cA(i)=\Set{T',T''}$ for distinct $T',T''$, where $T''\neq T$ but possibly $T'=T$.\COMMENT{suppose $T'\in \cS\cap \cA(i)$ such that $T''\in (\cS-\Set{T})\cap \cA(i)\cap \cT_{T'}(i)$. By definition of $\cT_{T'}(i)$, there is $\cS'\in \cX_{T',j',j'-4}(i)$, where $4\le j'\le j_{max}$, such that $(\cS'-\Set{T'})\cap \cA(i)=\Set{T''}$.} Consider such a pair $\cS,\cS'$. In particular, $|V(\cS)\cap V(\cS')|\ge 4$. Let $H:=(\cS\cup \cS')\cap \cC(i)$. Since $H$ is obtained from $(\cS\cup \cS')-\Set{T,T''}$ by deleting $j-3-c$ edges,\COMMENT{$\cS'$ contains no available triples which are not already contained in $\cS$. $\cS$ has $j-2-(c-1)$ available triples, two of which are already deleted.} we deduce from Proposition~\ref{prop:butterfly count} and Fact~\ref{fact:extension induction}(i) that $(H,T)$ is $(j-3-c)$-balanced, and so $(H,T) \in Ext( \kappa , \ell )$ for some $\kappa\le j-3-c$ and $\ell \in [m -1]$.\COMMENT{$H$ is non-empty as $|\cS'\cap \cC(i)|\ge 2$}
%
%
%
By~\eqref{extension bound} and \eqref{extension error def}, we conclude that the number of such pairs is $\O(n^{j-3-c+\frac{m-1}{m}})$.
From~\eqref{erdos regularity}, we have that $|\cX_{T,j,c-1}(i)|=f_{j,c-1}(i) \pm \eps(i)n^{j-2-c}$. Using Lemma~\ref{lem:properties before stopping}\ref{properties before stopping number available} and~\ref{properties before stopping number threats}, we conclude that
\begin{eqnarray*}
E^{gain} &=&  \sum_{\cS\in \cX_{T,j,c-1}(i)} \frac{|((\cS-\Set{T})\cap \cA(i))\sm \bigcup_{T'\in \cS\cap \cA(i)}\cT_{T'}(i)|}{|\cA(i)\sm (\cT_T(i)\cup \Set{T})|} \\
         &=& \left(|\cX_{T,j,c-1}(i)|-\O(n^{j-2-c-\frac{1}{m}}) \right) \frac{j-2-c}{|\cA(i)\sm (\cT_T(i)\cup \Set{T})|}\\
				&=& (f_{j,c-1}(i) \pm 2\eps(i)n^{j-2-c})\frac{j-2-c}{A(i)\pm 2\eps(i) n^3}\\
				&\overset{\eqref{O fractions},\eqref{triple trajectory magnitude}}{=}& \frac{j-2-c}{A(i)}f_{j,c-1}(i) + \O_{\gamma}(\eps(i))n^{j-5-c}.
\end{eqnarray*}

Thus, using Lemma~\ref{lem:derivatives}\ref{lem:derivatives erdos}, \ref{lem:derivatives erdos zero} and~\ref{erdos trajectory taylor}, we obtain
\begin{align}
\cexpn{\tilde{L}}{\Delta X_{T,j,c}(i) \mid i<\tau_{freeze,T}}  &= -E^{loss} + E^{gain} =f'_{j,c}(i) + \O_{\gamma}(\eps(i))n^{j-5-c} \nonumber\\
                        &=\Delta f_{j,c}(i) + \O_{\gamma}(\eps(i))n^{j-5-c}. \label{expected erdos trend}
\end{align}

We infer that
\begin{eqnarray*}
\cexpn{\tilde{L}}{\Delta X^{\pm}_{T,j,c}(i) \mid i<\tau_{freeze,T}} &\overset{\eqref{margin def formal triple}}{=}& \pm \cexpn{\tilde{L}}{\Delta X_{T,j,c}(i) \mid i<\tau_{freeze,T}} \\
                                                     & &\mp \Delta f_{j,c}(i) -\Delta \eps(i) n^{j-3-c} \\
																&\overset{\eqref{expected erdos trend},\eqref{error function props}}{=}&  \O_{\gamma}(\eps(i))n^{j-5-c} - C \eps(i) n^{j-5-c} \overset{\eqref{hierarchy}}{\le} 0.											
\end{eqnarray*}
With the observation at \eqref{supermartingale sufficient}, this shows that $(X_{T,j,c}^\pm(0),X_{T,j,c}^\pm(1),\dots)$ is a supermartingale with respect to $(L(0),L(1),\dots)$.

Moreover, using \eqref{triple trajectory magnitude}, \eqref{edge and danger magnitude} and~\eqref{constant density cutoff}, we can also deduce that 
\begin{align}
\cexpn{\tilde{L}}{|\Delta X_{T,j,c}(i)| \mid i<\tau_{freeze,T}} \le  E^{loss} + E^{gain} = \O_{\gamma}(n^{j-5-c}).\label{abs bound erdos}
\end{align}

\begin{step}
The expected change of $X_{R,(H,U)}$
\end{step}

Finally, consider $\ell\in \Set{1,\dots,m}$, $\kappa \in\Set{0,\dots,\ell}$, $(H,U)\in Ext(\kappa,\ell)$ and $R\In V$ with $|R|=|U|$. By Fact~\ref{fact:extension induction}(i), we have that $\kappa(H-e,U)\le \kappa+1$ for all $e\in H$. Thus, using~\eqref{extension bound}, Lemma~\ref{lem:properties before stopping}\ref{properties before stopping number available} and the fact that $i<\tau_{stop}$, we obtain
\begin{align}
\cexpn{\tilde{L}}{\Delta X_{R,(H,U)}(i)} &\le \sum_{e\in H} \frac{X_{R,(H-e,U)}(i)}{|\cA(i)|} \le \frac{|H| 2n^{\kappa+1+\frac{\ell}{m+\kappa+1}}}{A(i)-\eps(i)n^3}  \nonumber \\
                                           &=    \O_\gamma(n^{\kappa-2+\frac{\ell}{m+\kappa+1}}) \le n^{\kappa-2+\frac{\ell}{m+\kappa}}. \label{expected extension trend}
\end{align}
(Here, we use Fact~\ref{fact:extension type trivial}\ref{fact:extension type empty} and~\ref{fact:extension type trivial bound} instead of \eqref{extension bound} if $H-e$ is empty.)

We continue to obtain
\begin{eqnarray*}
\cexpn{\tilde{L}}{\Delta X^+_{R,(H,U)}(i) } &\overset{\eqref{margin def formal extension}}{=}& \cexpn{\tilde{L}}{\Delta X_{R,(H,U)}(i) } - \Delta \eps_{\kappa,\ell}(i)  \\
                                           & \overset{\eqref{expected extension trend},\eqref{extension error def}}{\le} & n^{\kappa-2+\frac{\ell}{m+\kappa}} -  n^{\kappa-2+\frac{\ell}{m+\kappa}} =0.
\end{eqnarray*}
By the observation at~\eqref{supermartingale sufficient}, $(X_{R,(H,U)}^+(0),X_{R,(H,U)}^+(1),\dots)$ is a supermartingale with respect to $(L(0),L(1),\dots)$.

Moreover, since $\Delta X_{R,(H,U)}(i)\ge 0$, we immediately have that 
\begin{align}
\cexpn{\tilde{L}}{|\Delta X_{R,(H,U)}(i)|  } \le n^{\kappa-2+\frac{\ell}{m+\kappa}}. \label{abs bound extension}
\end{align}

\begin{step}
Expected absolute changes
\end{step}

From \eqref{abs bound edge}, \eqref{abs bound erdos} and~\eqref{abs bound extension}, it is now easy to deduce with the triangle inequality, Lemma~\ref{lem:derivatives}\ref{edge trajectory taylor},\ref{erdos trajectory taylor} and \eqref{error function props}, \eqref{extension error def} that
\begin{align*}
\cexpn{\tilde{L}}{|\Delta X^{\pm}_e(i)| \mid i<\tau_{freeze,e}} &= \O_\gamma(n^{-1}), \\
\cexpn{\tilde{L}}{|\Delta X^\pm_{T,j,c}(i)| \mid i<\tau_{freeze,T}} &=\O_{\gamma}(n^{j-5-c}), \\
\cexpn{\tilde{L}}{|\Delta X^+_{R,(H,U)}(i)| }  & \le 2n^{\kappa-2+\frac{\ell}{m+\kappa}}.
\end{align*}
With~\eqref{expected abs sufficient}, this completes the proof.
\endproof

\subsection{Boundedness hypotheses} \label{subsec:boundedness}

We now establish boundedness hypotheses for the variables we track.

\begin{lemma} \label{lem:boundedness edge}
For every $2$-set $e$, we have $\Delta X^\pm_e(i)=\O(n^{\frac{1}{2}})$ for all $i$.
\end{lemma}

\proof
Fix a $2$-set $e$. For $i\ge \tau_{freeze,e}$, we trivially have $\Delta X^\pm_e(i)=0$. Suppose that $i< \tau_{freeze,e}$. In particular, $e\in E(i)$ and $e\in E(i+1)$. Note that $\cX_e(i+1)\sm \cX_e(i)=\emptyset$ and $\cX_e(i)\sm \cX_e(i+1)=\set{T\in \cX_e(i)}{T \leftrightarrow T^\ast(i)}$. Thus, $$|\Delta X_e(i)|\le \max_{T^\ast \in \cA(i)\sm \cX_e(i)}|\set{T\in \cX_e(i)}{T \leftrightarrow T^\ast}|.$$
Fix any $T^\ast \in \cA(i)\sm \cX_e(i)$. It follows that $|T^\ast \cup e|\ge 4$. The number of $T\in \cX_e(i)$ with $T \leftrightarrow T^\ast$ is bounded from above by the number of \Erd-configurations $\cS$ on $j\le j_{max}$ points with $e\cup T^\ast \In V(\cS)$, $|\cS\cap \cC(i)|=j-4$ and $|(\cS-\cC(i))[e\cup T^\ast]|\ge 1$, which by Corollary~\ref{cor:Erdos extension applicable} is $\O(n^{\frac{1}{2}})$.%
\COMMENT{If $j=4$ then we cannot apply Corollary~\ref{cor:Erdos extension applicable}, but in this case $V(\cS)=e\cup T^*$, so there are only $\O(1)$ choices.}

It follows that $\Delta X_e(i)=\O(n^{\frac{1}{2}})$, which, via \eqref{margin def formal edge}, implies $\Delta X^\pm_e(i)=\O(n^{\frac{1}{2}})$ using Lemma~\ref{lem:derivatives}\ref{edge trajectory taylor} and~\eqref{error function props}.
\endproof

\begin{lemma} \label{lem:boundedness triple}
For every $3$-set $T$, all $j\in \Set{6,\dots,j_{max}}$ and $c\in \Set{0,\dots,j-4}$, we have $\Delta X^\pm_{T,j,c}(i)=\O(n^{j-3-c-\frac{1}{m}})$ for all $i$.
\end{lemma}

\proof
For $i\ge \tau_{freeze,T}$, we trivially have $\Delta X^\pm_{T,j,c}(i)=0$. Suppose that $i< \tau_{freeze,T}$. In particular, $T\in \cA(i)$ and $T\in \cA(i+1)$.
We first examine the maximum positive change of $X_{T,j,c}(i)$. Note that if $c=0$, we clearly have $\Delta X_{T,j,c}(i)\le 0$. If $c>0$, we have
\begin{align*}
\Delta X_{T,j,c}(i) &\le \max_{T^\ast \in \cA(i)\sm \Set{T}}|\set{\cS\in \cX_{T,j,c-1}(i)}{T^\ast\in \cS}|\\
                 &\le \max_{T^\ast \in \cA(i)\sm \Set{T}}|\set{\cS\in \mathfrak{J}_j}{T,T^\ast\in \cS, |\cS\cap \cC(i)|=c-1}|\\
								&=\O(n^{j-3-(c-1)-2+\frac{1}{2}})=\O(n^{j-4-c+\frac{1}{2}})
\end{align*}
by Corollary~\ref{cor:Erdos extension applicable} (with $T\cup T^\ast$, $2$, $c-1$ playing the roles of $R,a,c$).

We next examine the maximum negative change of $X_{T,j,c}(i)$. Note that
\begin{align*}
-\Delta X_{T,j,c}(i) &\le \max_{T^\ast \in \cA(i)\sm \Set{T}}|\set{\cS\in \cX_{T,j,c}(i)}{T^\ast\mbox{ threatens }\cS,T}|.
\end{align*}
Fix $T^\ast\in \cA(i)\sm \Set{T}$. Suppose first that $c<j-4$. Observe that $|\cT_{T^\ast}(i)|=\O(n)$ by Lemma~\ref{lem:properties before stopping}\ref{properties before stopping number threats}. Thus, using Corollary~\ref{cor:Erdos extension applicable}, we obtain
\begin{align*}
|\set{\cS\in \cX_{T,j,c}(i)}{T^\ast\mbox{ threatens }\cS,T}| &\le \sum_{T'\in (\cT_{T^\ast}\cup \Set{T^\ast})\sm\Set{T}}|\set{\cS\in \mathfrak{J}_j}{T,T'\in \cS, |\cS\cap \cC(i)|=c}|\\
                       &=\O(n)\cdot \O(n^{\max\Set{j-3-c-2,0} +\frac{1}{2}}) = \O(n^{j-4-c+\frac{1}{2}}),
\end{align*}
as desired.
If $c=j-4$, then we have $|\set{\cS\in \cX_{T,j,j-4}(i)}{T^\ast\mbox{ threatens }\cS,T}|\le X_{T,T^\ast}(i) \le \O(n^{1-\frac{1}{m}})$ by Fact~\ref{fact:negative danger bound} and Corollary~\ref{cor:double extension}.

We conclude that $\Delta X_{T,j,c}(i)=\O(n^{j-3-c-\frac{1}{m}})$, which, via \eqref{margin def formal triple}, implies $\Delta X^\pm_{T,j,c}(i)=\O(n^{j-3-c-\frac{1}{m}})$ using Lemma~\ref{lem:derivatives}\ref{erdos trajectory taylor} and~\eqref{error function props}.
\endproof

\begin{lemma} \label{lem:boundedness extension}
Let $\ell\in \Set{1,\dots,m}$, $\kappa \in\Set{0,\dots,\ell}$, $(H,U)\in Ext(\kappa,\ell)$ and $R\In V$ with $|R|=|U|$. Then $\Delta X^+_{R,(H,U)}(i)=\O(n^{\kappa+\frac{\ell-1}{m+\kappa}})$ for all $i$.
\end{lemma}

\proof
For $i\ge \tau_{stop}$, we trivially have $\Delta X^+_{R,(H,U)}(i)=0$. Suppose that $i< \tau_{stop}$. We bound the maximum change of $X_{R,(H,U)}(i)$.

By Fact~\ref{fact:extension induction}(ii), we have that $\kappa(H-H[U\cup e],U \cup e)\le \kappa$ for all $e\in H$.
Fix any $T^\ast \in \cA(i)$. We need to give an upper bound on the number of $\phi\colon V(H)\to V$ which are $(H,U)$-extensions at $R$ in $\cC(i)\cup \Set{T^\ast}$, but not in $\cC(i)$. Fix any such $\phi$. Then there must be $e\in H$ with $\phi(e)=T^\ast$, and $\phi(e')\in \cC(i)$ for all $e'\in H-\Set{e}$. Thus, we have that $\phi$ is an $(H-H[U\cup e],U\cup e)$-extension at $R\cup T^\ast$ in $\cC(i)$.\COMMENT{$\phi(U\cup e)=\phi(U)\cup \phi(e)=R\cup T^\ast$ and for all $e'\in H-H[U\cup e]$, $\phi(e')\in \cC(i)$} The number $|V(H)\sm (U\cup e)|$ of free vertices in the new extension type is at most $\ell-1$ since $e\not\In U$.

Hence, by~\eqref{extension bound}, the number of all possible $\phi$ is at most $$\sum_{e\in H\colon |U \cup e|=|R\cup T^\ast|} X_{R\cup T^\ast,(H-H[U\cup e],U\cup e)}(i)= \O(n^{\kappa+\frac{\ell-1}{m+\kappa}}).$$
(Here, we use Fact~\ref{fact:extension type trivial}\ref{fact:extension type empty} and~\ref{fact:extension type trivial bound} instead of \eqref{extension bound} if $H-H[U\cup e]$ is empty.)
This implies $\Delta X^+_{R,(H,U)}(i)=\O(n^{\kappa+\frac{\ell-1}{m+\kappa}})$ since $\Delta \eps_{\kappa,\ell}(i)=n^{\kappa-2+\frac{\ell}{m+\kappa}}$ by~\eqref{extension error def}.
\endproof

\subsection{Proof of Lemma~\ref{lem:negative trend variables}}

We now prove Lemma~\ref{lem:negative trend variables}, which in turn implies Theorem~\ref{thm:process stop} and hence Theorem~\ref{thm:approx STS}.

\lateproof{Lemma~\ref{lem:negative trend variables}}
Fix a $2$-set $e$. By Lemma~\ref{lem:supermartingales}, $X^\pm_e$ form supermartingales, and $\expn{|\Delta X^\pm_e(i)| \mid  L(i)}=\O_{\gamma}(n^{-1})$ for all $i$. By Lemma~\ref{lem:boundedness edge} we have $\Delta X^\pm_e(i)=\O(n^{\frac{1}{2}})$ for all $i$. 
By~\eqref{initial edge}, we have $-X^\pm_e(0)=\Omega_{\eps_0}(n)$. Thus, we can apply~\eqref{freedman compact} with $(\alpha_1,\alpha_2,\alpha_3)=(1,\frac{1}{2},-1)$ to conclude that $\prob{\exists i\colon X^\pm_e(i)\ge 0}\le \eul^{-\Omega_{\eps_0}(n^{1/2})}$.

Fix a triple $T$, $j\in\Set{6,\dots,j_{max}}$ and $c\in \Set{0,\dots,j-4}$. By Lemma~\ref{lem:supermartingales}, $X^\pm_{T,j,c}$ form supermartingales, and $\expn{|\Delta X^\pm_{T,j,c}(i)| \mid  L(i)}=\O_{\gamma}(n^{j-5-c})$ for all $i$. By Lemma~\ref{lem:boundedness triple} we have $\Delta X^\pm_{T,j,c}(i)=\O(n^{j-3-c-\frac{1}{m}})$ for all $i$. By~\eqref{initial triple}, we have $-X^\pm_{T,j,c}(0)=\Omega_{\eps_0}(n^{j-3-c})$. Thus, we can apply~\eqref{freedman compact} with $(\alpha_1,\alpha_2,\alpha_3)=(j-3-c,j-3-c-\frac{1}{m},j-5-c)$ to conclude that $\prob{\exists i\colon X^\pm_{T,j,c}(i)\ge 0}\le \eul^{-\Omega_{\eps_0}(n^{1/m})}$.

Fix $\ell\in \Set{1,\dots,m}$, $\kappa \in\Set{0,\dots,\ell}$, $(H,U)\in Ext(\kappa,\ell)$ and $R\In V$ with $|R|=|U|$. By Lemma~\ref{lem:supermartingales}, $X^+_{R,(H,U)}$ forms a supermartingale, and $\expn{|\Delta X^+_{R,(H,U)}(i)| \mid  L(i)}= \O(n^{\kappa-2+\frac{\ell}{m+\kappa}})$ for all $i$. By Lemma~\ref{lem:boundedness extension}, we have $\Delta X^+_{R,(H,U)}(i)=\O(n^{\kappa+\frac{\ell-1}{m+\kappa}})$ for all $i$. Since $-X^+_{R,(H,U)}(0)=n^{\kappa+\frac{\ell}{m+\kappa}}$ by~\eqref{initial extension}, we can apply~\eqref{freedman compact} with $(\alpha_1,\alpha_2,\alpha_3)=(\kappa+\frac{\ell}{m+\kappa},\kappa+\frac{\ell-1}{m+\kappa},\kappa-2+\frac{\ell}{m+\kappa})$ to conclude that $\prob{\exists i\colon X^+_{R,(H,U)}\ge 0}\le \eul^{-\Omega(n^{1/2m})}$.\COMMENT{$\alpha_1-\alpha_2=\frac{1}{m+\kappa}\ge 1/2m$}

Thus, a final union bound shows that whp, all the variables $X^\pm$ are non-positive.
\endproof

\section{Counting sparse Steiner triple systems} \label{sec:counting}

Wilson conjectured that the number $STS(n)$ of non-isomorphic Steiner triple systems on $n$ vertices (provided $n$ is admissible) is $(n/\eul^2+o(n))^{n^2/6}$. This was recently proved by Keevash~\cite{keevash:18}.
Letting $STS_k(n)$ denote the number of $k$-sparse Steiner triple systems on $n$ vertices, we expect from our heuristics in Section~\ref{subsec:heuristics} that 
\begin{align}
	\log{STS_k(n)} \approx \log{STS(n)} - \int_{0}^{n^2/6} \rho(i),
\end{align}
where $\rho(i)=\sum_{j=6}^{k+2}\frac{J_j}{\binom{n}{3}^{j-3}}i^{j-3}$.
Let $erd_j$ denote the number of unlabelled \Erd-configurations on $[j]$ containing the triple $123$. Thus, $J_j=erd_j\binom{n-3}{j-3}$.

We thus have 
\begin{align}
	\int_{0}^{n^2/6} \rho(i) \approx \frac{n^2}{6}\sum_{j=6}^{k+2} \frac{erd_j}{(j-2)!}.
\end{align}

Hence, we conjecture that
\begin{align}
	STS_k(n) &= \left(n \eul^{-2-\sum_{j=6}^{k+2} \frac{erd_j}{(j-2)!}}+o(n) \right)^{\frac{n^2}{6}} .
\end{align}
In particular, since $erd_6=6$, we conjecture that the number of Pasch-free Steiner triple systems is $\left(n \eul^{-2-1/4}+o(n) \right)^{\frac{n^2}{6}} $. It would be interesting to find out whether the upper bound could be established using the entropy method as in \cite{LL:13}.

\section{General sparse designs} \label{sec:designs}

In this section, we discuss the possible existence of sparse Steiner systems with more general parameters.
Given $n\ge q>r\ge 2$, a \defn{partial $(n,q,r)$-Steiner system} is a set $\cS$ of $q$-subsets of some $n$-set $V$ such that every $r$-subset of $V$ is contained in at most one $q$-set in~$\cS$. 
An \defn{$(n,q,r)$-Steiner system} is a partial $(n,q,r)$-Steiner system $\cS$ with $|\cS|=\binom{n}{r}/\binom{q}{r}$, i.e.~every $r$-set is covered.
For fixed $q$ and $r$, we call $n$ \defn{admissible} if $\binom{q-i}{r-i}\mid \binom{n-i}{r-i}$ for all $0\le i\le r-1$. It is easy to see that this condition is necessary for the existence of an $(n,q,r)$-Steiner system. Recently, Keevash~\cite{keevash:14} was able to settle the so-called existence conjecture, stating that for sufficiently large $n$, there exists an $(n,q,r)$-Steiner system whenever $n$ is admissible (see~\cite{GKLO:16} for an alternative proof).

\subsection{A generalized Erd\H{o}s-conjecture}
In order to formulate a generalized Erd\H{o}s-conjecture, we first consider what the generalized \Erd-configurations might be.
A \defn{$(j,\ell)_{q,r}$-configuration} is a set of $\ell$ $q$-sets on $j$ points every two of which intersect in at most $r-1$ points.

The reason why $(j,j-3)$-configurations appear in every Steiner triple system $\cS$ 
is that whenever we have a $(j,\ell)$-configuration $\cL$ in $\cS$ with an uncovered pair, then the triple in $\cS$ which covers this pair determines only one new point, i.e.~we can extend $\cL$ to a $(j+1,\ell+1)$-configuration. Since there trivially are $(4,1)$-configurations to start with, we can obtain $(j,j-3)$-configurations for all $j\ge 4$ in this way. For more general Steiner systems, the argument is similar. Having already found some configuration, we consider an uncovered $r$-set inside this configuration, and the Steiner system returns a $q$-set which covers this $r$-set. Hence, we add one $q$-set on the expense of maximally $q-r$ new points. 
This leads to the following family of functions: $$\kappa_{q,r}(j):=\lflr \frac{j-r-1}{q-r} \rflr.$$ 
Note that $\kappa_{3,2}(j)=j-3$, and more generally, $\kappa_{r+1,r}(j)=j-r-1$.

\begin{prop}\label{prop:generalized erdos}
For all $n\ge j > q>r\ge 2$, every $(n,q,r)$-Steiner system $\cS$ contains a $(j,\kappa_{q,r}(j))_{q,r}$-configuration.
\end{prop}

\proof
We first prove by induction on $x\in \bN_0$ that the statement holds for all $j$ of the form $j=x(q-r)+q+1$ (for which we have $\kappa_{q,r}(j)=x+1$). For $x=0$, we have $j=q+1$ and $\kappa_{q,r}(j)=1$, and can thus just take any $q$-set of $\cS$ together with an arbitrary additional point.

Suppose now that $x\ge 1$. Let $j':=(x-1)(q-r)+q+1$. Note that $\kappa_{q,r}(j')=x$. By induction, there is a $(j',x)_{q,r}$-configuration $\cL$ in $\cS$. Since $$x\binom{q}{r}<\binom{(x-1)(q-r)+q+1}{r}= \binom{j'}{r},$$ there exists an $r$-set $e\In V(\cL)$ with $e\notin \cL$.\COMMENT{Need that $$xq(q-1)\cdots(q-r+1)< ((x-1)(q-r)+q+1)((x-1)(q-r)+q)\cdots ((x-1)(q-r)+q-r+2).$$ It's enough to have $x(q-r+2)(q-r+1)< ((x-1)(q-r)+q-r+3)((x-1)(q-r)+q-r+2)= (x(q-r)+3)(x(q-r)+2)$, which is equivalent to $(x^2-x)(q-r)^2+2x(q-r-1)+6>0$.}

This is covered by a unique $q$-set $Q$ of $\cS$. Thus, $\cL\cup \Set{Q}$ is a collection of $x+1$ $q$-sets of $\cS$ on at most $j'+q-r$ points.\COMMENT{That every two of them intersect in at most $r-1$ points follows because they're all in $\cS$.} 
By adding isolated vertices if necessary, we may assume that this yields a $(j'+q-r,x+1)_{q,r}$-configuration, i.e.~a $(j,\kappa_{q,r}(j))_{q,r}$-configuration.

For general $j$, write $j=x(q-r)+q+1+y$ for $x\in\bN_0$ and $0\le y<q-r$. Then $\kappa_{q,r}(j)=x+1=\kappa_{q,r}(j-y)$. Thus, by the above, there exists a $(j-y,\kappa_{q,r}(j))_{q,r}$-configuration, and we may simply add isolated vertices to obtain a $(j,\kappa_{q,r}(j))_{q,r}$-configuration.
\endproof

Proposition~\ref{prop:generalized erdos} tells us that we cannot forbid $(j,\kappa_{q,r}(j))_{q,r}$-configurations.
Motivated by this, we say that an $(n,q,r)$-Steiner system $\cS$ is \defn{$k$-sparse} if there is no $(j,\kappa_{q,r}(j)+1)_{q,r}$-configuration in $\cS$ with $2\le \kappa_{q,r}(j)+1\le k$. Note that this coincides with the definition of $k$-sparseness for triple systems. We propose the following generalization of \Erd's conjecture.

\begin{conj} \label{conj:new Erdos}
For all $q>r\ge 2$ and every $k$, there exists an $n_k$ such that for all admissible $n>n_k$, there exists a $k$-sparse $(n,q,r)$-Steiner system.
\end{conj}
For the case $r=2$, this has already been conjectured in~\cite{FR:13}.\COMMENT{Conjecture 1.4}

\subsection{Partial result}

It is not clear why the proof of Proposition~\ref{prop:generalized erdos} would yield the `correct' function $\kappa_{q,r}$. We now provide some evidence that $\kappa_{q,r}$ is indeed the correct function. It would be interesting to see whether our process can be generalized to prove Conjecture~\ref{conj:new Erdos} approximately. We take a much simpler route here and show that if we allow even one more $q$-set per $j$ vertices, then the conjecture approximately holds.
We say that a (partial) $(n,q,r)$-Steiner system $\cS$ is \defn{weakly $k$-sparse} if there is no $(j,\kappa_{q,r}(j)+2)_{q,r}$-configuration in $\cS$ with $\kappa_{q,r}(j)+2\le k$.

As tools, we use the Lov\'asz local lemma and a result on almost perfect matchings in hypergraphs due to Pippenger. The idea is to first randomly sparsify the set of $q$-sets in such a way that no $j$-set contains too many $q$-sets, whilst preserving certain degree and codegree conditions. This allows to find an almost perfect matching in a suitable auxiliary hypergraph, producing an approximate Steiner system which is automatically sparse.

For events $B_1,\dots,B_n$ in a common probability space, we say that the graph $\Gamma$ with $V(\Gamma)=[n]$ is a \defn{dependency graph} if $B_i$ is mutually independent of all $B_j$ with $ij\notin \Gamma$.

\begin{lemma}[Lov\'asz local lemma, cf.~\cite{AS:08}]\label{lem:LLL}
Let $B_1,\dots,B_n$ be events with dependency graph~$\Gamma$. If there exist $x_1,\dots,x_n \in [0,1)$ such that for all $i\in [n]$, we have
\begin{align}
\prob{B_i}\le x_i \prod_{j\in N_{\Gamma}(i)}(1-x_j), \label{LLL condition}
\end{align}
then $$\prob{\bigcap_{i=1}^n\overline{B_i}}\ge \prod_{i=1}^n(1-x_{i}).$$
\end{lemma}

The following is a well-known result due to Pippenger, which has never been published, but several stronger versions have been proven since (see e.g.~\cite{PS:89}).
\begin{theorem}[Pippenger\COMMENT{As stated in Alon-Kim-Spencer}] \label{thm:PP nibble}
Suppose $1/D, \eps \ll \gamma,1/r$. Let $H$ be an $r$-graph on $n$ vertices and suppose that
$d_H(x)=(1 \pm \eps)D$ for all $x\in V(H)$ and $d_H(\Set{x,y})\le \eps D$ for all distinct $x,y\in V(H)$. Then there exists a matching in $H$ covering all but $\gamma n$ vertices.
\end{theorem}

\begin{theorem}
Let $1/n\ll \gamma,1/k,1/q$ and $2\le r<q$. There exists a weakly $k$-sparse partial $(n,q,r)$-Steiner system $\cS$ on $n$ vertices with $|\cS|\ge (1- \gamma) \binom{n}{r}/\binom{q}{r}$.
\end{theorem}

\proof
Let $V$ be a set of size $n$. Choose new constants $\eps,\theta$ such that $1/n\ll \eps,\theta \ll \gamma,1/k,1/q$.
Let $j_{max}$ be the maximal $j$ such that $\kappa_{q,r}(j)+2\le k$. Thus, we may assume that $\frac{j-r}{q-r}+1 \ge \frac{j-r+\theta}{q-r-\theta}$ for all $q+1\le j\le j_{max}$.
Note that $\kappa_{q,r}(j)\ge \frac{j-r}{q-r}-1$ and thus we have
\begin{align}
(q-r-\theta)(\kappa_{q,r}(j)+2)\ge j-r+\theta \label{theta choice}
\end{align}
for all $q+1\le j\le j_{max}$.
Let $\cA$ be the random $q$-graph on $V$ obtained by selecting every $Q\in \binom{V}{q}$ independently with probability $p:=n^{-(q-r)+\theta}$.

For $q+1\le j\le j_{max}$ and a set $S\in \binom{V}{j}$, we let $B_S$ denote the event that $|\cA[S]|\ge \kappa_{q,r}(j)+2$. 
For an $r$-set $e\In V$, let $B_e$ be the event that $d_{\cA}(e)\neq (1\pm \eps)n^\theta/(q-r)! $. 
Finally, let $V_{codeg}$ be the set of all pairs $e,e'$ of distinct $r$-sets in $V$. For $ee'\in V_{codeg}$, let $B_{ee'}$ be the event that $d_{\cA}(e\cup e')\ge n^{\theta/10} $.

We claim that with positive probability, none of the events $B_S$, $B_e$, $B_{ee'}$ occurs. (This will allow us to apply Theorem~\ref{thm:PP nibble}.)
Define the graph $\Gamma$ with vertex set $V(\Gamma)=\bigcup_{j=q+1}^{j_{max}}\binom{V}{j}\cup \binom{V}{r}\cup V_{codeg}$ and add the following edges: add a clique on $\binom{V}{r}\cup V_{codeg}$, and for $S,S'\in \bigcup_{j=q+1}^{j_{max}}\binom{V}{j}$, $e\in \binom{V}{r}$ and $e_1e_2\in V_{codeg}$, add
\begin{align}
SS'\in E(\Gamma) & \mbox{ if }|S\cap S'|\ge q, \nonumber\\
eS\in E(\Gamma) &\mbox{ if } e\In S, \nonumber\\
\Set{e_1,e_2}S\in E(\Gamma) &\mbox{ if } e_1\cup e_2\In S.\nonumber
\end{align}
Clearly, $\Gamma$ is a dependency graph for the events $(B_v)_{v\in V(\Gamma)}$. We now aim to fulfill the conditions of the Lov\'asz local lemma.

Clearly, for $q+1\le j\le j_{max}$ and $S\in \binom{V}{j}$, we have
\begin{align}
\prob{B_S}\le \O(1) p^{\kappa_{q,r}(j)+2}. \label{prob bad S}
\end{align}
Now, consider $e\in \binom{V}{r}$.
Note that $\expn{d_{\cA}(e)}=p\binom{n-r}{q-r}$. Using a standard Chernoff-Hoeffding bound,\COMMENT{
Let $X$ be the sum of independent Bernoulli random variables. Then the following hold.
\begin{enumerate}[label={\rm(\roman*)}]
\item For all $0\le\eps \le 3/2$, $\prob{X \neq (1\pm \eps)\expn{X} } \leq 2\eul^{-\eps^2\expn{X}/3}$.
\item If $t\ge 7 \expn{X}$, then $\prob{X\ge t}\le \eul^{-t}$.
\end{enumerate}
Use (i) here (ii) for the next one} we have that
\begin{align}
\prob{B_e} \le \eul^{-\frac{\eps^2}{4}n^{\theta}}. \label{prob bad deg}
\end{align}
Finally, consider $ee'\in V_{codeg}$. Note that $\expn{d_{\cA}(e\cup e')}\le p n^{q-r-1}=o(1)$. Thus, using a standard Chernoff-Hoeffding bound, we have that
\begin{align}
\prob{B_{ee'}} \le \eul^{-n^{\theta/10}}. \label{prob bad codeg}
\end{align}

For all $q+1\le j\le j_{max}$ and $S\in \binom{V}{j}$, define $x_S:=x_j:=n^{-j+r}$. For all $e\in \binom{V}{r}$, define $x_e:=x_{deg}:=\eul^{-\frac{\eps^2}{4}n^{\theta/2}}$. For all $ee'\in V_{codeg}$, define $x_{ee'}:=x_{codeg}:=\eul^{-n^{\theta/20}}$.

We now check condition~\eqref{LLL condition}. First consider $q+1\le j\le j_{max}$ and $S\in \binom{V}{j}$. We have
\begin{align*}
\prod_{v\in N_{\Gamma}(S)} (1-x_v) \ge \prod_{j'=q+1}^{j_{max}}(1-x_{j'})^{\binom{j}{q}n^{j'-q}} \cdot (1-x_{deg})^{\binom{j}{r}} \cdot (1-x_{codeg})^{\binom{j}{r}^2} \ge 1/2.
\end{align*}
\COMMENT{each term tends to $1$}
Since
\begin{align*}
	\prob{B_S} \overset{\eqref{prob bad S}}{\le} \O(1) p^{\kappa_{q,r}(j)+2} = \O(1) n^{-(q-r-\theta)(\kappa_{q,r}(j)+2)} \overset{\eqref{theta choice}}{\le} \frac{n^{-(j-r)}}{2},
\end{align*}
we conclude that \eqref{LLL condition} is satisfied for $S$.
Now consider $e\in \binom{V}{r}$. We have
\begin{align*}
\prod_{v\in N_{\Gamma}(e)} (1-x_v) \ge \prod_{j=q+1}^{j_{max}}(1-x_{j})^{n^{j-r}} \cdot (1-x_{deg})^{n^r} \cdot (1-x_{codeg})^{n^{2r}} \ge \eul^{-j_{max}}.
\end{align*}
\COMMENT{$(1-x_{deg})^{n^r} \cdot (1-x_{codeg})^{n^{2r}}\to 1 $, $(1-x_{j})^{n^{j-r}}\to \eul^{-1}$, so get $0.9\eul^{-(j_{max}-q)}$, say.}
Since $\prob{B_e}\le \eul^{-\frac{\eps^2}{4}n^{\theta}} \le \eul^{-j_{max}}x_{deg}$ by~\eqref{prob bad deg}, we deduce that \eqref{LLL condition} is satisfied for $e$.
A similar calculation also shows that  \eqref{LLL condition} is satisfied for $ee'\in V_{codeg}$.%
\COMMENT{We have
\begin{align*}
\prod_{v\in N_{\Gamma}(ee')} (1-x_v) \ge \prod_{j=q+1}^{j_{max}}(1-x_{j})^{n^{j-(r+1)}} \cdot (1-x_{deg})^{n^r} \cdot (1-x_{codeg})^{n^{2r}} \ge 0.9.
\end{align*}
Since $\prob{B_{ee'}}\le \eul^{-n^{\theta/10}} \le 0.9 x_{codeg}$ by~\eqref{prob bad codeg}, we deduce that \eqref{LLL condition} is satisfied for $ee'$.}

Thus, with Lemma~\ref{lem:LLL} we can infer that with positive probability none of the events $(B_v)_{v\in V(\Gamma)}$ occurs. Let $\cA$ be such a $q$-graph. Define the auxiliary $\binom{q}{r}$-graph $H$ with $V(H)=\binom{V}{r}$ and $E(H)=\set{\binom{Q}{r}}{Q\in \cA}$. Since no event $B_e$ occurred, we have $d_{H}(e)=(1\pm \eps)n^\theta/(q-r)!$ for all $e\in V(H)$. Moreover, since no event $B_{ee'}$ occurred, we have $d_{H}(\Set{e,e'})\le n^{\theta/10}$ for all distinct $e,e'\in V(H)$.
Hence, by Theorem~\ref{thm:PP nibble}, there exists a matching in $H$ covering all but $\gamma \binom{n}{r}$ vertices of $H$. This correponds to a partial $(n,q,r)$-Steiner system $\cS$ on $V$ covering all but $\gamma \binom{n}{r}$ $r$-sets. Thus, $|\cS|\ge (1- \gamma) \binom{n}{r}/\binom{q}{r}$. Finally, every $q$-set of $\cS$ is contained in $\cA$, and since no event $B_S$ occurred, $\cS$ is weakly $k$-sparse.
\endproof

\section*{Acknowledgement}

We are grateful to Tom Bohman and Lutz Warnke for pointing out a minor oversight in the calculation of $E^{gain}$ in an earlier version of this paper.


\vspace{1cm}

{\footnotesize \obeylines \parindent=0pt

Stefan Glock, Daniela K\"{u}hn, Allan Lo, Deryk Osthus
\vspace{0.3cm}
School of Mathematics
University of Birmingham
Edgbaston
Birmingham
B15 2TT
UK
}
\vspace{0.3cm}
\begin{flushleft}
{\it{E-mail addresses}:}
\tt{[s.glock,d.kuhn,s.a.lo,d.osthus]@bham.ac.uk}
\end{flushleft}

\end{document}